%% file: nccg1.tex
\documentclass[twoside,12pt,a4paper]{article}
\usepackage[marked,reporte]{myart98}
\ppnum{funct-an/9703001}{1997}{}
\input{mydef}

\newcommand{\SL}{SL(2,\Space{R}{})}
\newcommand{\matr}[4]{{\ensuremath{ \left( \begin{array}{cc}
#1 & #2 \\ #3 & #4 \end{array}\right) }}}

\newcommand{\oper}[1]{\ensuremath{\mathcal{#1}}}

\titleshort{Two Approaches to Non-Commutative Geometry}
\authorshort{Vladimir V. Kisil}

\begin{document}
\title{Two Approaches to \\Non-Commutative Geometry\thanks{Supported
by grant 3GP03196 of the FWO-Vlaanderen (Fund of Scientific 
Research-Flanders),
Scientific Research Network ``Fundamental Methods and Technique in
Mathematics''.}}
\author{Vladimir V. Kisil\\
		Institute of Mathematics\\
               Economics and Mechanics\\
               Odessa State University\\
               ul. Petra Velikogo, 2\\
               Odessa-57, 270057, UKRAINE\\
E-mail: \texttt{kisilv@member.ams.org}}
\date{April 1, 1997; Revised February 3, 1998}
\maketitle
\begin{abstract}
Looking to the history of mathematics one could find out two 
outer approaches to Geometry. First one (algebraic) is due to 
\person{Descartes} 
and second one (group-theoretic)---to \person{Klein}. We will see 
that they are not rivalling but are tied (by \person{Galois}). We also 
examine their modern life as philosophies of non-commutative 
geometry. Connections between different objects (see \emph{keywords})
are discussed.
\AMSMSC{46H30}{30G35, 47A13, 81R05}
\end{abstract}
\newpage
\begin{small}
\tableofcontents
\end{small}  
\section{Introduction}
\epigraph{Descartes was the first great modern philosopher, a founder 
of modern biology, a first-rate physicist, and only incidentally a 
mathematician.}{M.~Kline \cite[\S~15.3]{Kline72}.}{}
If one takes a look on mathematics as a collection of facts about a 
large diversity of objects covered by the Mathematic Subject 
Classification by digits from $00$ till $93$, then it will be 
difficult to explain why we are referring to mathematics as a ``united 
and inseparable'' science. Even 
the common origin of all mathematical fields laid somewhere in 
\emph{Elements} of \person{Euclid} hardly be an excuse. Indeed western 
philosophy and psychology as well commonly rooted in works of Great 
Greeks (\person{Plato} or \person{Aristotle}?), but it will be offensive
for a philosopher as well as for a psychologist do not 
distinguish them nowadays. In contrast, in 
mathematics the most exiting and welcome results usually combine 
facts, notions and ideas of very remote fields (the already classical 
example is \person{Atiah}-\person{Singer} theorem linking algebra, 
analysis and geometry).

So even while exponentially growing mathematical facts remain linked by
some underlying ideas, which are essentially vital for the mathematical
enterprise. Thus it seems worth enough for a mathematician to think (and
sometime to write) not only about \emph{new}
mathematical facts but also about \emph{old} original ideas.
Particularly relations between new facts and old ideas deserve 
special attention.

In this paper we return to two great ideas familiar to every 
mathematician: \emph{Coordinate method} of \person{R.~Descartes} and 
\emph{Erlangen program} of \person{F.~Klein}. They specify two 
different approaches to geometry and become actual now in connections 
with  such an interesting and fascinate area of research as 
non-commutative geometry~\cite{Connes94,Manin91}.

The paper format is as follows. After some preliminary remarks in 
Sections~\ref{se:coord}--\ref{se:manin} we will arrive to paper's core 
in Section~\ref{se:analytic}: it is an abstract scheme of  
analytic function theories. We will interpret classic 
examples accordingly to the presented scheme in Section~\ref{se:names}. 
In Section~\ref{se:calculus} we will expand consideration to 
non-commutative spaces with physical application presented in
Section~\ref{se:observ}. For a reader primary interested in analytic 
function theory Section~\ref{se:analytic} and~\ref{se:names} 
form an independent reading.

A feature of the paper is a healthy portion of self-irony. It will not 
be an abuse if a reader will smile during the reading as the author did
during the writing.

\section{Coordinates: from Descartes to Nowadays} \label{se:coord}
\epigraph{\ldots{} coordinate geometry changed the face of 
mathematics.}{M.~Kline \cite[\S~15.6]{Kline72}.}{}
The first danger of a deep split in mathematics condensed before 17-th 
century. The split shown up between synthetic geometry of Greeks and 
abstract algebra of Arabs. It 
was overcome by works of \person{Fermat} and \person{Descartes} on 
analytic geometry. Descartes pointed out that geometric construction 
call for adding, subtracting, multiplying, and dividing segments of 
lines. But this is exactly four operations of algebra. So one could
express and solve geometric problems in algebraic terms.
The crucial steps are: 
\begin{enumerate}

\item\label{it:coord} To introduce a coordinate system, which allows 
us to label geometrical points by sets of numbers or 
\emph{coordinates}. 
Of course, coordinates are \emph{functions} of points.

\item\label{it:equations} To link desired geometrical properties and 
problems with \emph{equations} of coordinate. Then solutions to 
equations will give solutions of original geometric problems.
\end{enumerate}
Four algebraic operations suggest to consider coordinates from 
Step~\ref{it:coord} as members of a commutative algebra of functions 
of geometric points. 
Moreover, it seems unavoidable to use other functions from this 
algebra to construct and solve equations from Step~\ref{it:equations}.
In such a way the commutative algebra of functions become the central 
object of new analytic geometry. From geometry this algebra inherits 
some additional structures (metrical, topological, etc.)

Having a vocabulary that allows us to translate from the geometrical 
language to algebraic one nothing could prevent one to use it in the 
opposite direction: \emph{what are geometrical counterparts of such and 
such algebraic notions}? Particularly, what will happen if we start 
from an arbitrary commutative algebra (with a relevant additional 
structure) instead coordinate one? 
The answer to a question of this type is given by the celebrated theorem
of \person{I.M.~Gelfand} and 
\person{M.A.~Naimark}~\cite[\S~4.3]{Kirillov76}:
\begin{thm}[Gelfand-Naimark]
Any commutative Banach algebra $\algebra{A}$ 
could be realized as the 
algebra $\FSpace{C}{}(X)$ of all continuous functions over a 
topological space $X$. The points of $X$ are labelled by the maximal 
ideals (or characters) of $\algebra{A}$.
\end{thm}

So looking for adventures and new experience in geometry one should 
try to begin from \emph{non-commutative} algebras. This gives the 
name of resulting theory as \emph{non-commutative geometry}.
Very attractive (and 
in some sense archetypical)
object is \emph{anticommutative} coordinates, which satisfy to identity 
$x_1 x_2 =- x_2 x_1$ instead of commutation rule $x_1 x_2 =x_2 x_1$:
\begin{quotation}
In recent years, several attempts have been made to develop 
anticommutative analogs of geometric notions. Most of these attempts are 
rooted in one of the great ideas of mathematics in this century, the 
idea that geometry is recovered from algebra by taking prime or 
maximal ideals of a commutative ring. Thus one replaces a commutative 
ring by an exterior algebra and one tries to do something 
similar.~\cite{ChanRotaStein95} 
\end{quotation}
On this road one could found natural unification of vector geometry 
and Grassmann algebras.

There are still many other opportunities. As it always happens with 
generalizations of a non-trivial notion, there is no \emph{the} 
non-commutative geometry but there \emph{are} non-commutative 
geometries (for example, \cite{Connes94} and~\cite{Manin91}). 
\begin{rem}
It is important that all described constructions deserve the name of 
\emph{geometry} not only by peculiar rules of the game named 
mathematics but also by their relations with the real word. Their 
ambitions are to describe the structure of actual space(-time) on the 
quantum level (see Section~\ref{se:observ}) exactly as classic geometry 
succeed it for the Nile valley. Thus as commutative geometry is 
the language of classic mechanics so non-commutative geometry is 
supposed to be a language for quantum physics. So in modern science 
non-commutative and quantum are usually synonymous.
\end{rem}
\begin{example}\label{ex:quantum}
Quantum mechanics begun from the Heisenberg commutation 
relations $XY=YX +i\hbar I$ and was a source of inspiration 
for functional analysis. Thus an additional structure on non-commutative
algebra, which describes geometrical objects, comes often from operator
theory. 
Table~\ref{ta:c-q} contains a vocabulary for non-commutative geometry
taken from~\cite{ChamConnes96b}, which is devoted to description of 
quantized theory of gravitational field.\\ 
\begin{table}[t]
\begin{center}
\begin{tabular}{ccc}
\textbf{Commutative - Classic
} &\qquad &\textbf{Non-Commutative - Quantum}\\
Complex variables & & Operator in Hilbert space\\
Real variables & & Selfadjoint operators\\
Infinitesimal variable & & Compact operator\\
Infinitesimal of order $\alpha$ & & $\mu_n(T)=O(n^{-\alpha})$\\
Integral & &$\int T=\log$ trace $T$
\end{tabular}
\end{center}
\caption{Vocabulary between Commutative-Classic and 
Non-Commutative-Quantum Geometries. Here $\mu_n(T)$ is stand for
$n$-th eigenvalue of a compact operator $T$.} \label{ta:c-q}
\end{table}
On this way one could find such exiting examples as a 
curvature or a connection on a geometric set consisting of two points 
only~\cite{Connes94}. We will discuss such a vocabulary with more 
details in Section~\ref{se:observ}.
\end{example}

\section{Klein, Galois, Descartes}\label{se:klein}
\epigraph{\ldots{} not all of geometry can be incorporated in 
Klein's scheme.}{M.~Kline \cite[\S~38.5]{Kline72}.}{}
The second (group-theoretic) viewpoint on geometry was expressed by 
\person{F.~Klein} in his famous \emph{Erlangen program} and seems 
unrelated to coordinate geometry. But this is not true.
\begin{thm} The difference between Klein and Descartes is Galois or
\begin{equation}\label{eq:kgd}
Klein(1872)=Galois(1831)+ Descartes(1637)
\end{equation}
\end{thm}
\begin{proof}
The idea of Descartes to reduce geometric problems to the algebraic 
equations was a brilliant prophecy regarding to the very modest 
knowledge about algebraic equations in that times. Although his 
contemporaries could solve already algebraic equations up to fourth 
degree, i.e., the most general that could be solved in radicals, the 
general theory of algebraic equation was not even supposed to exist. 
It is enough to say that for Descartes\footnote{We could be surprised 
how little a genius should be acquainted with mathematical facts to 
discover the underlying fundamental idea.} even negative roots of 
an equation were false roots not speaking about complex 
ones~\cite[\S13.5]{Kline72}.

We were lucky that another genius, namely \person{Galois}, leaves to us 
a message at the end of his very short life: \emph{solvability of 
equations could be determined via its group of symmetries}. Again we 
could skip all facts of corresponding theory (for elementary account 
see for example~\cite[\S~31.4]{Kline72}) and just add it to the 
coordinate geometry\footnote{Note the appearance of a geometric 
problem in~\cite[\S~31.5]{Kline72} immediately after the Galois theory 
in~\cite[\S~31.4]{Kline72}}. The result will almost identical with the 
Erlangen program~\cite[\S~38.5]{Kline72}:

\begin{defn} \label{de:geometry}
Every geometry can be characterized by a \emph{group} (of 
transformations) and that geometry is really concerned with 
\emph{invariants} under this group.
\end{defn}

To see this achievement in the proper light one should keep in mind 
that the notion of abstract group was fully recognized only years 
after Erlangen program~\cite[\S~49.2]{Kline72} (as the notion of a root
to an equation was recognized years after Descartes' coordinate
geometry).
\end{proof}
\begin{cor}
Klein is greater than Descartes:
\begin{displaymath}
Klein > Descartes
\end{displaymath}
\end{cor}
\begin{proof}
Galois is obviously positive ($Galois>0$), so the~\eqref{eq:kgd} 
implies the assertion.
\end{proof}
\begin{prob}
Find a counterexample to the last Corollary.
\end{prob}
\emph{Hint}. See the epigraph.

Even not being universal the Erlangen program offers a way to 
classify the variety of geometries by means of associated
groups~\cite[\S~38.5]{Kline72}. We
will link in the Erlangen spirit different analytic function theories 
and connected functional calculi with group representations. This gives 
an opportunity to make their classification too.
\begin{rem}
The connection jokingly presented in this section is often 
unobserved but is very important. In concrete situations it takes  
practical forms. Many connections between non-commutative geometry and
group representations are mentioned in~\cite{Segal96b}. For example, an
early result on non-commutative geometry---the classification of
von~Neumann factors of type III---was carried out by means their group
of automorphisms~\cite{Connes94} (see also Remarks~\ref{re:disk}
and~\ref{re:calc}). And that is also important, such a connection
generates a hope that not only particular mathematical facts are linked by
fundamental ideas, but also these ideas are connected one with another on 
a deeper level.
\end{rem}

\section{The Heisenberg Group and Symmetries of the Manin Plain}
\label{se:manin}
\epigraph{\ldots{} I want to popularize the Heisenberg group, which 
is remarkably little known considering its 
ubiquity.}{R.~Howe \cite{Howe80a}.}{}

What does Erlangen program tell to a non-commutative geometer? 
One could do non-commutative geometry considering (quantum) symmetries 
of non-commutative (=quantum) objects. This approach leads to quantum 
groups~\cite{Kassel94,Manin88a,Manin91}. We will consider a simple 
example of $M_q(2)$.

We have mentioned already (see Remark~\ref{ex:quantum}) the Heisenberg 
commutation relation:
\begin{equation}\label{eq:a-heisenberg}
[X,Y]=I.
\end{equation}
Here $X$, $Y$, and $I$, viewed as generators of a Lie algebra, 
produce a Lie group, which is the celebrated Heisenberg (or Weyl) 
group $\Heisen{1}$~\cite{Folland89,Howe80a,MTaylor86}. An element 
$g\in \Heisen{n}$ (for any positive integer $n$) could be represented 
as $g=(t,z)$ with 
$t\in\Space{R}{}$, $z=(z_1,\ldots,z_n)\in \Space{C}{n}$ and the group 
law is given by  
\begin{equation}\label{eq:g-heisenberg}
g*g'=(t,z)*(t',z')=(t+t'+\frac{1}{2}\sum_{j=1}^n\Im(\bar{z}_j 
z_j'), z+z'),
\end{equation}
where $\Im z$ denotes the imaginary part of a complex number $z$.
Of course the Heisenberg group is non-commutative and particularly one 
could find out for $\Heisen{1}$:
\begin{equation}\label{eq:weyl}
(0,1)*(0,i)=(1,0)*(0,i)*(0,1)
\end{equation}
The last formula is usually called the Weyl commutation relation and 
is the integrated (or exponentiated) form of~\eqref{eq:a-heisenberg}.

Now let us consider a continuous irreducible 
representation~\cite{Kirillov76,MTaylor86}   
$\pi:\Heisen{1}\rightarrow\algebra{A}$ of the Heisenberg group in a 
\Cstar-algebra $\algebra{A}$. Then elements of the center of 
$\Heisen{1}$, which 
are of the form $(t,0)$, will map to multipliers $ae$ of the identity 
$e\in\algebra{A}$. For unitary representations we  have 
$\modulus{a}=1$ for all such $a$.

Let $x=\pi(0,1)$, $y=\pi(0,i)$, $qe=\pi(1,0)$. 
Then the Weyl commutation relation~\eqref{eq:weyl} will correspond to
\begin{equation}\label{eq:manin}
xy=qyx
\end{equation}
The last equation is taken by \person{Yu. Manin}~\cite{Manin88a} as 
the defining relation for \emph{quantum plane}, which is also often 
called the Manin plane. As we saw the Manin plane is a 
representation of the Heisenberg-Weyl group and thus could be 
named\footnote{Now and then the nobility of a mathematical object is 
measured by the number of its family names.}  
Heisenberg-Weyl-Manin plane. More formal: Manin plane with a parameter 
$q\in\Space{C}{}$ is the quotient of free algebra generated by elements 
$x$ and $y$ subject to two-sided ideal generated by the quadratic 
relation~\eqref{eq:manin}.

Let us consider the regular representation $\pi_r$ of 
$\Heisen{1}$ as right shifts on $\FSpace{L}{_2}(\Heisen{1})$ and let 
$\algebra{M}$ be an algebra of operators generated by 
$x_r=\pi_r(0,1)$ and $y_r=\pi_r(0,i)$. Then Manin plane is a 
representation $\pi$ of $\algebra{M}$ under which $x=\pi(x_r)$ and 
$y=\pi(y_r)$.  Moreover \emph{an algebraic identity $f(x,y)=0$ holds 
on Manin plane for all $q\in\Space{C}{}$ if and only if the algebraic 
identity $f(x_r,y_r)=0$ is true in $\algebra{M}$.}

Connection between the Manin plane and the Weyl commutation relation 
was already mentioned in~\cite{Manin88a} but was not used explicitly. 
Thus this connection disappeared from the following works. Even 
fundamental treatise~\cite{Kassel94} (which has 531 pages!) mentioned 
the Heisenberg group only twice and in both cases disconnected with 
the Manin plane. This is especially pity because the Heisenberg group 
could give new incites for quantum groups, as it does for the Manin 
plane. 

\begin{example}{}\footnote{For brevity the Example is presented very 
informal. We believe that this is not an abuse.} 
As it was mentioned at the beginning of the section the quantum group 
could be viewed as symmetries of quantum object. Particularly $M_q(2)$ 
is a set of $2\times 2$ matrixes $\matr{a}{b}{c}{d}$,
``which maps the Manin plane to the Manin plane'', namely if $x$ and 
$y$ satisfy to~\eqref{eq:manin} and
\begin{displaymath}
\matr{a}{b}{c}{d} \left(\begin{array}{c}
x\\y
\end{array}\right)=
\left(\begin{array}{c}
ax+by\\cx+dy
\end{array}\right)
\end{displaymath}
then $x'=ax+by$, $y'=cx+dy$ should again satisfy to~\eqref{eq:manin}.
From this condition for both matrixes  $\matr{a}{b}{c}{d}$ and  
$\matr{a}{c}{b}{d}$ one could find that entries $a$, $b$, $c$,
and $d$ of $\matr{a}{b}{c}{d}$ are subject to the following six 
identities:
\begin{eqnarray}
ab&=&q^{-1}ba,\quad ac=q^{-1}ca,\quad db=qbd,\quad dc=qcd \label{eq:6-1}\\
bc&=&cb,\quad [a,d]=-(q-q^{-1})cb \label{eq:6-2}
\end{eqnarray}
To keep anytime in mind all six identities for any element of 
$M_q(2)$ is a little bit laborious, isn't it? The Heisenberg group 
proposes an elegant way out. $2$-vector $\left(\begin{array}{c}
x\\y
\end{array}\right)$ on the Manin plane is of course an element of a 
representation of $\Heisen{2}$ as the Manin plane itself is a 
representation of $\Heisen{1}$. And it is easy to check that due to 
celebrated Stone-von Neumann theorem~\cite{Folland89,MTaylor86} six 
identities \eqref{eq:6-1}--\eqref{eq:6-2} up to unitary equivalence 
prescribe that $a$, $b$, $c$, $d$  belong to a representation $\pi_M$ of
$\Space{H}{2} $ such that:
\begin{displaymath}
a=\pi_M (0,1,1), \quad b=\pi_M (0,1,i), \quad c=\pi_M (0,i,1), \quad 
d=\pi_M (0,i,i).
\end{displaymath}
Then a product of two elements of $M_q(2)$ is naturally described via 
representation of $ \Space{H}{4} $ (because entries of different 
elements pairwise commutes) and so far.
Thus \emph{$M_q(2)$ is a representation of 
$M_2(\pi_M(\Heisen{\infty}))$}. There is no problem to think over $ 
\Space{H}{ \infty } $ because consideration is purely algebraic and no 
topology is involved. Alternatively one could consider $ \Space{H}{N}$ 
for some unspecified sufficiently large $N$.
Of course, our conclusion about algebraic identities on the Manin plane 
remains true for $M_q(2)$.
\end{example}

\section{Analytic Function Theory from Group Representations} 
\label{se:analytic}
\epigraph{I learned from I.M.~Gel'fand that Mathematics of any kind
is a representation theory\ldots}{Yu. I. Manin,\\}{Short address on 
occasion of Professor Gel'fand 70th birthday.}

One could notice that through glasses of Erlangen program the boundary 
between commutative and non-commutative geometry becomes tiny. 
Quantum groups became a representation of classic group and the word 
\emph{non-commutative} is not distinguishing any more. Indeed, the 
group of symmetries of Euclidean geometry is non-commutative (think on 
a composition of a shift and a rotation) and a group could 
act as symmetries of classic and quantum object (symplectic 
transformations). And even more: the difference between geometry 
and, say, analysis became unnoticeable.

In this section we are going to proof the following\footnote{We are
continuing the definition-theorem-proof style presentation of  
Section~\ref{se:klein}.} 
\begin{thm} Complex analysis $=$ conformal geometry.
\end{thm}
\begin{proof}
The \emph{group} of conformal mapping (particularly for the unit disk 
they are fraction-linear transformations) plays an important role in any 
textbook on complex analysis~\cite[Chap.~10]{Krantz82}, 
\cite[Chap.~2]{Rudin80}. Moreover,
two key object of complex analysis, namely the Cauchy-Riemann equations
and Cauchy kernel, are \emph{invariants} with respect to these
transformations.
Thus the Theorem follows from Definition~\ref{de:geometry}.
\end{proof}

\subsection{Wavelet Transform and Coherent States}

We agree with a reader if he/she is not satisfied by the last short proof
and would like to see a more detailed account how the core of complex
analysis could be reconstructed from representation theory of $\SL$.  We
present an abstract scheme, which also could be applied to other analytic
function theories~\cite{CnopsKisil97a,Kisil97c}.  We start from a dry
construction followed in the next Section by classic examples, which will
justify our usage of personal names.

Let $ X $ be a topological space and let $G$ be a group that acts $G: 
{X}\rightarrow {X}$ as a transformation $g: x \mapsto g \cdot x$ from
the left, i.e., $g_1 \cdot(g_2 \cdot x)=(g_1 g_2)\cdot 
x$. Moreover, let $G$ act on $X$ transitively. Let there
exist a measure $dx$ on $ X$ such that a representation $\pi(g): f(x)
\mapsto m(g,x) f( g^{-1}\cdot x)$ (with a function $m(g,h)$) is unitary
with respect to the scalar product
$\scalar{f_1(x)}{f_2(x)}_{\FSpace{L}{2}(X ) } =
\int_{  X} f_1(x) \bar{f}_2(x) \,d(x)$, i.e.,
\begin{displaymath} 
\scalar{[\pi(g)f_1](x)}{[\pi(g)f_2](x)}_{\FSpace{L}{2}( X ) }
= \scalar{f_1(x)}{f_2(x)}_{\FSpace{L}{2}( X ) }\qquad \forall f_1, 
f_2 \in \FSpace{L}{2}(X).
\end{displaymath}
We consider the Hilbert space $\FSpace{L}{2}( X )$ where representation
$\pi(g)$ acts by unitary operators.
\begin{rem}
It is well known that the most 
developed part of representation theory consider unitary 
representations in Hilbert spaces. By this reason we restrict our 
attention to Hilbert spaces of analytic functions, the area usually 
done by means of the functional analysis technique. We also assume that 
our functions are complex valued and this is sufficient for examples 
explicitly considered in the present paper. However the presented scheme
is working also for vector valued functions and this is the natural 
environment for Clifford analysis~\cite{BraDelSom82}, for example. One 
also could start from an abstract Hilbert space $H$ with no explicit 
realization as $ \FSpace{L}{2}(X)$ given.
\end{rem}  

Let $H$ be a closed compact\footnote{While the compactness will be 
explicitly used during our abstract consideration, it is not crucial 
in fact. Example~\ref{ex:segal} will show how one could make a trick for 
non-compact $H$.} subgroup of $G$ and let $f_0(x)$ be such a
function that $H$ acts on it as the multiplication
\begin{equation} \label{eq:homogenious}
[\pi(h)f_0](x)=\chi(h) f_0(x) \qquad , \forall h\in H,
\end{equation}
by a function $\chi(h)$, which is a character of $H$
i.e., $f_0(0)$ is a common eigenfunction for all operators
$\pi(h)$. Equivalently $f_0(x)$ is a common eigenfunction for operators 
corresponding under $\pi$ to a basis of the Lie algebra of $H$.
Note also that $ \modulus{\chi(h)}^2=1 $ because $\pi$ is unitary.
$f_0(x)$ is called
\emph{vacuum vector} (with respect to subgroup $H$). 
We introduce the 
$\FSpace{F}{2}(X)$ to be the closed liner subspace of 
$\FSpace{L}{2}(X)$ uniquely defined 
by the conditions:
\begin{enumerate}
\item\label{it:begin} $f_0\in \FSpace{F}{2}( X )$;
\item $\FSpace{F}{2}( X )$ is $G$-invariant;
\item\label{it:end} $\FSpace{F}{2}( X )$ is 
$G$-irreducible.
\end{enumerate}
Thus restriction of $\pi$ on $ \FSpace{F}{2}(X) $ is an irreducible 
unitary representation.

The \emph{wavelet transform}\footnote{The subject of coherent states or
wavelets have been arising many times in many \emph{applied} areas and
the author is not able to give a comprehensive history and proper
credits. One could mention important
books~\cite{Daubechies92,KlaSkag85,Perelomov86}. We give our
references by recent paper~\cite{Kisil95a}, where applications to
\emph{pure} mathematics were considered.} $\oper{W}$ could be defined 
for square-integral representations $\pi$ by the formula~\cite{Kisil95a}
\begin{eqnarray}
\oper{W}&:& \FSpace{F}{2}( X ) \rightarrow
\FSpace{L}{ \infty }(G) \nonumber \\
&:& f(x) \mapsto
\widetilde{f}(g)=\scalar{f(x)}{\pi(g)f_0(x)}_{\FSpace{L}{2}(
X ) } \label{eq:wavelets}
\end{eqnarray}
The principal advantage of the wavelet transform $\oper{W}$ is that it 
express the representation $\pi$ in geometrical terms. Namely it 
\emph{intertwins} $\pi$ and left regular representation $ \lambda $ on 
$G$:
\begin{equation} \label{eq:g-inter}
[\lambda_g\oper{W} f](g')
=[\oper{W}f](g^{-1}g')
= \scalar{f}{\pi_{g^{-1}g'}f_0}
= \scalar{\pi_g f}{\pi_{g'}f_0}
=[\oper{W}\pi_g f](g'),
\end{equation}
i.e., $\lambda \oper{W} = \oper{W} \pi$. Another important feature of 
$W$ is that it does not lose information, namely function $ f(x) $ could 
\emph{be recovered} as the linear combination of \emph{coherent states} 
$f_g(x)=[\pi_g f_0](x)$ from its wavelet transform $ \widetilde{f}(g)
$~\cite{Kisil95a}: 
\begin{equation} \label{eq:g-inverse}
f(x)=\int_G \widetilde{f}(g) f_g(x)\, dg=\int_G \widetilde{f}(g) [\pi_g f_0](x) 
\,dg,
\end{equation}
where $dg$ is the Haar measure on $G$ normalized such that 
$\int_G \modulus{ \widetilde{f}_0(g)}^2\,dg=1 $. One also has an
orthogonal \emph{projection} $ \widetilde{\oper{P}}$ from $
\FSpace{L}{2}(G,dg)$ to image $ \FSpace{F}{2}(G,dg)$ of $
\FSpace{F}{2}(X) $ under wavelet transform $\oper{W}$, which is just a
convolution on $g$ with the image $ \widetilde{f}_0(g)=\oper{W}(f_0(x))$
of the vacuum vector~\cite{Kisil95a}:
\begin{equation} \label{eq:g-proj}
[ \widetilde{ \oper{P}}w](g')=\int_{G} w(g)  \widetilde{f}_0(g^{-1}g')\,dg.
\end{equation}

\subsection{Reduced Wavelets Transform}
Our main observation will be that one could be much more economical (if
subgroup $H$ is non-trivial) with a help of~\eqref{eq:homogenious}: in
this case one need to know $ \widetilde{f}(g) $ not on the whole group $G$
but only on the homogeneous space $G/H$~\cite[\S~3]{AliAntGazMue}.

Let $ \Omega=G / H$ and $s: \Omega \rightarrow G$ be a continuous
mapping~\cite[\S~13.1]{Kirillov76}. Then any $g\in G$ has a unique
decomposition of the form $g=s(a)h$, $a\in \Omega$ and we will write
$a=s^{-1}(g)$, $h=r(g)={(s^{-1}(g))}^{-1}g$. Note that $ \Omega $ is a 
left $G$-homogeneous space\footnote{$ \Omega $ with binary operation $
(a_1,a_2) \mapsto s^{-1}(s(a_1)\cdot s(a_2))$ becomes a loop of the most 
general form~\cite{Sabinin72}. Thus theory of reduced wavelet transform 
developed in this subsection could be considered as \emph{wavelet 
transform associated with loops}. However we prefer to develop our 
theory based on groups rather on loops.}
with an action defined in terms of $s$ as follow: $g: a \mapsto
s^{-1}(g\cdot s(a)) $. Due to~\eqref{eq:homogenious} one could
rewrite~\eqref{eq:wavelets} as:
\begin{eqnarray*}
\widetilde{f}(g) & = & \scalar{f(x)}{\pi(g)f_0(x)}_{\FSpace{L}{2}( X ) }\\
& = & \scalar{f(x)}{\pi(s(a)h)f_0(x)}_{\FSpace{L}{2}( X ) }\\
& = & \scalar{f(x)}{\pi(s(a))\pi(h)f_0(x)}_{\FSpace{L}{2}( X ) }\\
& = & \scalar{f(x)}{\pi(s(a))\chi(h)f_0(x)}_{\FSpace{L}{2}( X ) }\\
& = & \bar{\chi}(h)\scalar{f(x)}{\pi(s(a))f_0(x)}_{\FSpace{L}{2}( X ) }
\end{eqnarray*}
Thus $\widetilde{f}(g)= \bar{\chi}(h)\widehat{f}(a)$ where 
\begin{equation} \label{eq:berg-cauch}
\widehat{f}(a)=[\oper{C} f] (a)=\scalar{f(x)}{\pi(s(a))f_0(x)}_{\FSpace{L}{2}(
X ) }
\end{equation}
and function $\widetilde{f}(g)$ on $G$ is completely defined by function $ 
\widehat{f}(a) $ on $ \Omega $. Formula~\eqref{eq:berg-cauch} gives us
an embedding $ \oper{C}: \FSpace{F}{2}(X) \rightarrow 
\FSpace{L}{ \infty }( \Omega ) $, which we will call \emph{reduced 
wavelet transform}. We denote by $ \FSpace{F}{2}(\Omega)$ the image of 
$ \oper{C}$ equipped with Hilbert space inner product induced 
by $ \oper{C}$ from $ \FSpace{F}{2}(X) $.

Note a special property of $ \widetilde{f}_0(g)$
and $ \widehat{f}_0(a) $:
\begin{displaymath}
\widetilde{f}_0(h^{-1}g)
= \scalar{f_0}{\pi_{h^{-1}g}f_0}
= \scalar{\pi_{h}f_0}{\pi_{g}f_0}
= \scalar{\chi(h)f_0}{\pi_{g}f_0}
=\chi(h) \widetilde{f}_0(g).
\end{displaymath}
It follows from~\eqref{eq:g-inter} that $ \oper{C} $ 
intertwines $ \rho\oper{C}= \oper{C} \pi$
representation $\pi$ with the representation
\begin{equation} \label{eq:a-inter}
[\rho_g \widehat{f}](a) = \widehat{f}(s^{-1}(g\cdot s(a))) 
\chi(r(g\cdot s(a))).
\end{equation}
While $\rho$ is not completely geometrical as $\lambda$ in applications 
it is still more geometrical than original $\pi$. In many cases 
$\rho$ is \emph{representat induced} by the character $\chi$.

If $ f_0(x) $ is a vacuum state with respect to $H$ then
$f_g(x)=\chi(h) f_{s(a)} (x)$ and  we could rewrite~\eqref{eq:g-inverse}
as follows:
\begin{eqnarray*}
f(x) & = &  \int_G \widetilde{f}(g) f_g(x)\, dg \nonumber \\
 & = & \int_{ \Omega } \int_H \widetilde{f}(s(a)h) f_{s(a)h}(x)\, 
 dh \, da \nonumber \\
 & = & \int_{ \Omega }  \int_H \widehat{f}(a) \bar{\chi}(h)  \chi(h) 
f_{s(a)} (x)   \, dh\,  da \nonumber \\ 
& = & \int_{ \Omega } \widehat{f}(a) f_{s(a)} (x)\,da \cdot
 \int_H \modulus{ \chi(h) }^2\, dh  \nonumber \\ 
 & = & \int_{ \Omega } \widehat{f}(a)  f_{s(a)}(x) \, da  ,
\end{eqnarray*}
if the Haar measure $dh$ on $H$ is set in such a way that 
$\int_H  \modulus{ \chi(h) }^2  \, dh=1$ and $dg=dh\,da$. We define
an integral transformation $ \oper{F} $ according to the last formula:
\begin{equation} 
[ \oper{F} \widehat{f}](x) =  \int_{ \Omega } \widehat{f}(a) 
f_{s(a)}(x) \, da \label{eq:a-inverse} ,
\end{equation}
which has the property $ \oper{F} \oper{C} = I$ on $ \FSpace{F}{2}(X) $
with $ \oper{C} $ defined in~\eqref{eq:berg-cauch}. One could consider
the integral transformation 
\begin{equation} \label{eq:szego}
[\oper{P} f](x)=[\oper{F} \oper{C} f](x)=
\int_{ \Omega } \scalar{f(y)}{f_{s(a)}(y)}_{\FSpace{L}{2}(
X ) }  f_{s(a)}(x) \, da
\end{equation}
as defined on whole $ \FSpace{L}{2}(X)$ (not only $ \FSpace{F}{2}(X) $).
It is known that $ \oper{P} $ \emph{is an orthogonal projection $
\FSpace{L}{2}(X) \rightarrow \FSpace{F}{2}(X) $}~\cite{Kisil95a}. 
If we formally use linearity of the scalar product $
\scalar{\cdot}{\cdot}_{ \FSpace{L}{2}(X)}$ (i.e., assume that the Fubini
theorem holds) we could obtain from~\eqref{eq:szego}
\begin{eqnarray}
[\oper{P} f](x) & = & 
\int_{ \Omega } \scalar{f(y)}{f_{s(a)}(y)}_{\FSpace{L}{2}(
X ) }  f_{s(a)}(x) \, da \nonumber \\
& = &  \scalar{f(y)}{\int_{ \Omega }f_{s(a)}(y) \bar{f}_{s(a)}(x) \, da
}_{\FSpace{L}{2}(
X )}   \nonumber \\
& = & \int_X f(y) K(y,x)\, d\mu(y) \label{eq:bergman},
\end{eqnarray}
where
\begin{displaymath}
K(y,x)=\int_{ \Omega } \bar{ f}_{s(a)}(y)  f_{s(a)}(x) \, da
\end{displaymath}
With the ``probability $ \frac{1}{2}$'' (see discussion on the Bergman and 
the Szeg\"o kernels bellow) the integral~\eqref{eq:bergman} exists in 
the standard sense, otherwise it is a singular integral operator (i.e, 
$K(y,x)$ is a regular function or a distribution).

Sometimes a reduced form $ \widehat{\oper{P}}: \FSpace{L}{2}( \Omega )
\rightarrow \FSpace{F}{2}( \Omega )$ of the 
projection $\widetilde{\oper{P}}$~\eqref{eq:g-proj} is of a separate
interest. It is an easy calculation that
\begin{eqnarray} \label{eq:s-b-proj}
[\widehat{\oper{P}} f](a')=\int_{ \Omega } f(a) \widehat{f}_0  ( 
s^{-1}( a^{-1}\cdot a')) \bar{\chi}(r(a^{-1}\cdot a')) \, da,
\end{eqnarray}
where $ a^{-1}\cdot a'$ is an informal abbreviation for $ {\left (s(a)
\right )}^{-1}\cdot s(a') $. As we will see its explicit form could be 
easily calculated in practical cases.

And only at the very end of our consideration we introduce the Taylor 
series and the Cauchy-Riemann equations. One knows that they are 
\emph{starting points} in the Weierstrass and the Cauchy 
approaches to complex analysis correspondingly.

For any decomposition $f_a(x)=\sum_\alpha  \psi_\alpha(x) V_\alpha(a)$ 
of the coherent states $f_a(x)$ by means of functions $V_\alpha(a)$
(where the sum could become eventually an integral) we have the
\emph{Taylor series} expansion
\begin{eqnarray} 
\widehat{f}(a) & = & \int_X f(x) \bar{f}_a(x)\, dx= \int_X f(x) \sum_\alpha 
\bar{\psi}_\alpha(x)\bar{V}_\alpha(a)\, dx  \nonumber \\
 & = &  \sum_\alpha 
\int_X f(x)\bar{\psi}_\alpha(x)\, dx \bar{V}_\alpha(a) \nonumber \\
 & = & \sum_{\alpha}^{\infty} \bar{V}_\alpha(a) f_\alpha,\label{eq:taylor}
\end{eqnarray}
where $f_\alpha=\int_X f(x)\bar{\psi}_\alpha(x)\, dx$.
However to be useful within the presented scheme such a decomposition 
should be connected with structures of $G$ and $H$. For example, if 
$G$ is a semisimple Lie group and $H$ its maximal compact subgroup then 
indices $\alpha$ run through the set of irreducible unitary
representations of $H$, which enter to the representation $\pi$ of $G$.

The \emph{Cauchy-Riemann equations} need more discussion. One could 
observe from~\eqref{eq:g-inter} that the image of $\oper{W}$ is
invariant under action of the left but right regular representations.
Thus $ \FSpace{F}{2}(\Omega)$ is invariant under 
representation~\eqref{eq:a-inter}, which is a pullback of the left 
regular representation on $G$, but its right counterpart. Thus generally
there is no way to define an action of left-invariant vector fields on 
$\Omega$, which are infinitesimal generators of right translations, on 
$ \FSpace{L}{2}(\Omega) $. But there is an exception. Let $ 
\algebra{X}_j $ be a maximal set of left-invariant vector fields on $G$
such that 
\begin{displaymath} 
\algebra{X}_j \widetilde{f}_0(g)=0.
\end{displaymath}
Because $\algebra{X}_j $ are left invariant we have $\algebra{X}_j 
\widetilde{f}_g'(g)=0$ for all $g'$ and thus image of $\oper{W}$, which 
the linear span of $\widetilde{f}_g'(g)$, belongs to intersection of 
kernels of $\algebra{X}_j$. The same remains true if we consider
pullback $ \widehat{\algebra{X}}_j$ of $\algebra{X}_j$ to $ \Omega$. 
Note that the number of linearly independent $\widehat{\algebra{X}}_j$ is 
generally less than for ${\algebra{X}}_j$. We call $ 
\widehat{\algebra{X}}_j$ as \emph{Cauchy-Riemann-Dirac} operators in 
connection with their property
\begin{equation} \label{eq:dirac}
\widehat{\algebra{X}}_j \widehat{f}(g)=0 \qquad \forall \widehat{f}(g) 
\in \FSpace{F}{2}(\Omega).
\end{equation}
Explicit constructions of the Dirac type operator for a discrete 
series representation could be found in 
\cite{AtiyahSchmid80,KnappWallach76}.

We do not use Cauchy-Riemann-Dirac operator in our construction, but 
this does not mean that it is useless. One could found at least such its 
nice properties:
\begin{enumerate}
\item Being a left-invariant operator it naturally encodes an
information about symmetry group $G$.
\item It effectively separates irreducible components of the 
representation $\pi$ of $G$ in $\FSpace{L}{2}(X )$.
\item It has a local nature in a neighborhood of a point vs. transformations, which act globally on the domain.
\end{enumerate}

\subsection{Great Names}
It is the time to give names of Greats to our formulas. A surprising thing 
is: all formulas have two names, each one for different 
circumstances (judgments will be given in Section~\ref{se:names}). 
More precisely we have two alternatives
\begin{table}[t]
\begin{center}
\begin{tabular}{|| c || c | c ||}  
\hline
\hline
Notion & Name for case  $ X \not\sim \Omega $ & Name for case $ 
X \sim \Omega $\\
\hline \hline
Space \FSpace{F}{2}(X) & Hardy space & 
\raisebox{-9pt }[0pt][0pt]{Bergman space} \\
\cline{1-2}
Space $\FSpace{F}{2}( \Omega )$ & Segal-Bargmann space & \\
\hline
Formula~\eqref{eq:berg-cauch} &
Cauchy Integral & \raisebox{-9pt }[0pt][0pt]{Bergman integral} \\
\cline{1-2}
Formula~\eqref{eq:a-inverse} & Segal-Bargmann inverse & \\
\hline
Projection~\eqref{eq:bergman} &
Szeg\"o projection & \raisebox{-9pt }[0pt][0pt]{Bergman projection}\\
\cline{1-2}
Projection~\eqref{eq:s-b-proj} &  Segal-Bargmann projection & \\
\hline
Series~\eqref{eq:taylor}& \multicolumn{2}{c||}{Taylor series}\\
\hline
Operator~\eqref{eq:dirac}& \multicolumn{2}{c||}{Cauchy-Riemann-Dirac 
operator} \\
\hline
\hline
\end{tabular}                  
\end{center}                  
\caption{Abstract formulas and their classical names. The apparent 
periodicity of names is not a coincidence or author's 
arbitrariness.} \label{ta:names}
\end{table}  
\begin{enumerate}
\item $ X $ and $ \Omega $ are \emph{not} isomorphic as topological 
spaces with measures. Then $ \FSpace{F}{2}(X)$ becomes the \emph{Hardy 
space}.
The space $ \FSpace{L}{2}( \Omega )$ is \emph{Segal-Bargmann space}.
We could consider~\eqref{eq:berg-cauch} as
an abstract analog of the \emph{Cauchy integral formula}, which recovers 
analytic function $ \widehat{f}(a) $ on $ \Omega $ from its (boundary) 
values $ f(x) $ on $ X $. The formula~\eqref{eq:a-inverse} could be 
named as \emph{Bargamann inverse formula}\footnote{The name is taken 
from the context of classical Segal-Bargmann
space~\cite{Bargmann61,Segal60} and stand for the mapping from the
Segal-Bargmann space to $ \FSpace{L}{2}( \Space{R}{n} )$. I do not know
its name in the context of Hardy space, but it definitely should exist
(it is impossible to find out  something new for analytic functions in 
our times).}. The orthogonal projection $ \oper{P} $ \eqref{eq:szego} is
the \emph{Szeg\"o projection} and expression~\eqref{eq:bergman} is a
singular integral operator. In this case of Hardy spaces we could not 
define a Hilbert space structure on $ \FSpace{F}{2}(\Omega) $ by 
means of a measure on $\Omega$. Therefor we need
to keep both faces of our theory: the scalar product could be defined by
a measure only for $ \FSpace{F}{2}(X)$ and analytic structure is defined
via geometry of $ \FSpace{F}{2}( \Omega )$.

\item $ X $ and $ \Omega $ \emph{are} isomorphic as topological spaces
with measures. In sharp contrast to the previous case where no name
appears twice, this case fully covered by single person, namely 
\person{Bergman}. On the other hand, because $ X \sim \Omega $ we have 
twice as less objects in this case. Particularly $ \FSpace{F}{2}(X) \sim 
\FSpace{F}{2}( \Omega )$ is known as the \emph{Bergman} space. In this
case representation $\pi$ on $ \FSpace{F}{2}(X) \sim \FSpace{F}{2}(
\Omega ) $ due to embedding $ s: \Omega \rightarrow G$ could be treated
as a part of \emph{left regular} representation on $ \FSpace{L}{2}(G) $
by shifts and thus representation $\pi$ belongs to \emph{discrete
series}~\cite[\S~VI.5]{Lang85}. The orthogonal projection $ \oper{P} $
\eqref{eq:szego} is the \emph{Bergman projection} and
expression~\eqref{eq:bergman} is a regular integral operator.

\end{enumerate}

\begin{rem}
The presented constructions being developed independently has many 
common points with the theory of harmonic functions on symmetric spaces
(see for example~\cite{Koranyi72a}). The new feature presented here is 
the association of function theories with sets of data $(G,\pi,H,f_0)$ 
rather symmetric spaces. This allows 
\begin{enumerate}
\item To distinguish under the common framework the Hardy and different
Bergman spaces, which are defined on the same symmetric space.
\item To deal with non-geometrical group action like in the 
case of the Segal-Bargmann space (see Example~\ref{ex:segal}).
\item To apply it not only to semisimple but also nilpotent Lie groups 
(see Example~\ref{ex:segal}).
\end{enumerate}
On the other hand \person{Bargmann} himself realizes that spaces of 
analytic functions give realizations of representations for semisimple
groups and later this idea was used many times (see for
example~\cite{AtiyahSchmid80}, \cite[\S~5.4]{GilbertMurray91}, and 
references therein). We would like emphasize here that not only analytic
functions spaces themselves come from representation theory. All
principal ingredients (the Cauchy integral formula, the Cauchy-Riemann 
equations, the Taylor series, etc.) of analytic function theory also 
appear from such an approach..
\end{rem}  

\section{The Name of the Game} \label{se:names}
\epigraph{It's the `ch' sound in Scottish words like \emph{loch} or 
German words like \emph{ach}; 
it's Spanish `j' and a Russian `kh'. When you say
it correctly to your computer, the terminal may become slightly
moist.}{Donald E. Knuth \cite[Chap.~1]{TeXbook}.}{}

It is the great time now to explain personal names appeared in the
previous Section. The main example is provided by group 
$G=\SL$ (books~\cite{HoweTan92,Lang85,MTaylor86} are our
standard references about $ \SL $ and its
representations) consisting of $ 2\times 2$ matrices $
\matr{a}{b}{c}{d}$ with real entries and determinant $ad-bc=1$. Via
identities
\begin{displaymath}
\alpha= \frac{1}{2}(a+d-ic+ib), \qquad \beta= \frac{1}{2}(c+b -ia+id)
\end{displaymath}
we have isomorphism of $ \SL $ with group $SU(1,1)$ of $
2\times 2$ matrices with complex entries of the form $ \matr{ \alpha}{
\beta}{ \bar{\beta}}{\bar{ \alpha}}$ such that $ \modulus{ \alpha }^2 - 
\modulus{ \beta }^2=1 $. We will use the last form for  $ \SL$ only.

$ \SL $ has the only non-trivial compact closed subgroup
$H$, namely the group of matrix of the form
$h_{\psi}=\matr{e^{i\psi}}{0}{0}{e^{-i\psi}}$.
Now any $ g\in \SL $ has a unique decomposition of the
form 
\begin{equation}
\matr{\alpha}{\beta}{\bar{\beta}}{\bar{\alpha}}=
\frac{1}{\sqrt{1- \modulus{a}^2 }}
\matr{1}{a}{\bar{a}}{1} 
\matr{e^{i\psi}}{0}{0}{e^{-i\psi}} 
\end{equation}
where $\psi=\Im \ln \alpha$, $a=\bar{\alpha}^{-1}\beta$, and $
\modulus{ a } < 1$ because $ \modulus{ \alpha }^2 - \modulus{ \beta
}^2=1 $. Thus we could identify $\SL / H$ with the unit
disk $ \Space{D}{} $ and define mapping $s: \Space{D}{} \rightarrow
\SL $ as follows
\begin{displaymath}
s: a \mapsto \frac{1}{\sqrt{1- \modulus{a}^2 }}
\matr{1}{a}{\bar{a}}{1}.
\end{displaymath}
The invariant measure $ d\mu(a) $ on $ \Space{D}{} $ coming from
decomposition $dg=d\mu(a)\, dh$, where $dg$ and $dh$ are Haar measures
on $G$ and $H$ respectively, is equal to 
\begin{equation}
d\mu(a)= \frac{da}{{(1- \modulus{a}^2)}^2}.
\end{equation}
with $da$---the standard Lebesgue measure on $ \Space{D}{} $.

The formula $g: a \mapsto g \cdot a= s^{-1}(g^{-1} s(a)) $ associates with 
a matrix $ g^{-1}=\matr{ \alpha}{ \beta}{ \bar{\beta}}{ \bar{\alpha}}$ the
fraction-linear transformation of $ \Space{D}{} $, which also could be 
considered as a transformation of $\dot{ \Space{C}{} }$ (the one-point 
compactification of $ \Space{C}{} $) of the form
\begin{equation} \label{eq:fr-lin}
g: z \mapsto g\cdot z=
\frac{\alpha z + \beta}{\bar{ \beta}z +\bar{\alpha}},  
\qquad g^{-1}=\matr{\alpha}{\beta}{\bar{\beta}}{\bar{\alpha}}.
\end{equation}
Complex plane $ \Space{C}{} $ is the disjoint union of three orbits,
which are acted by~\eqref{eq:fr-lin} of $G$ transitively:
\begin{enumerate}
\item The unit disk $ \Space{D}{}=\{ z \such \modulus{z}^2<1 \} $.
\item The unit circle $ \Gamma =\{ z \such \modulus{z}^2=1 \} $.
\item The remain part $ \Space{C}{} \setminus \bar{\Space{D}{}}
=\{ z \such \modulus{z}^2>1 \} $.
\end{enumerate}
$h_\psi \in H$ acts as the rotation on angle $2\psi$ and we could
identify $H \sim \Gamma $.

\begin{example}
First let us try $X=\Gamma$. We equip $ X $ with the standard Lebesgue
measure $ d\phi $ normalized in such a way that 
\begin{equation} \label{eq:n-measure}
\int_{ \Gamma } \modulus{f_0(\phi)}^2\, d\phi=1 \textrm{ with }
f_0(\phi)\equiv 1. 
\end{equation}
Then there is the unique (up to equivalence) way to make a unitary 
representation of  $G$ in $ \FSpace{L}{2}(X)= \FSpace{L}{2}( \Gamma, 
d\phi )$ from~\eqref{eq:fr-lin}, namely
\begin{equation} \label{eq:g-transform}
[\pi_g f](e^{i\phi})= \frac{1}{ \bar{\beta} e^{i\phi} + \bar{ \alpha }} 
f \left( \frac{  { \alpha }e^{i\phi}+{\beta}}{\bar{\beta} e^{i\phi} 
+ \bar{ \alpha }} \right).
\end{equation}
This is a realization of the \emph{mock discrete series} of $ \SL$.
Function $ f_0(e^{i\phi})\equiv 1 $ mentioned in~\eqref{eq:n-measure}
transforms as follows
\begin{equation} \label{eq:t-vacuum}
[\pi_g f_0](e^{i\phi})= \frac{1}{ \bar{\beta} e^{i\phi} + \bar{ \alpha }} 
\end{equation}
and particularly has an obvious property
$[\pi_{h_{\psi}}f_0](\phi)=e^{i\psi}
f_0(\phi)$, i.e., it is a \emph{vacuum vector} with respect to subgroup 
$H$. The smallest linear subspace $ \FSpace{F}{2}(X) \in
\FSpace{L}{2}(X) $ spanned by~\eqref{eq:t-vacuum} consists of boundaries
values of analytic function in the unit disk and is the \emph{Hardy}
space. Now reduced wavelet transform~\eqref{eq:berg-cauch} takes the
form
\begin{eqnarray}
\widehat{f}(a)=[\oper{C} f] (a) & =  &
\scalar{f(x)}{\pi(s(a))f_0(x)}_{\FSpace{L}{2}(
X ) } \nonumber \\
& = & \int_{\Gamma} f(e^{i\phi}) 
\frac{ \sqrt{ 1-\modulus{a}^2 } }{ \overline{ 
\bar{a}e^{i\phi} + 1} }\,d\phi \nonumber \\
& = &\frac{\sqrt{ 1-\modulus{a}^2 }}{i}  \int_{\Gamma} 
\frac{f(e^{i\phi})}{{a}+e^{i\phi}} ie^{i\phi}\,d\phi 
\nonumber \\
& = & \frac{\sqrt{ 1-\modulus{a}^2 }}{i}
 \int_{\Gamma}  \frac{f(z)}{{a}+z}\,dz, \label{eq:cauchy}
\end{eqnarray}
where $z=e^{i\phi}$.  Of course~\eqref{eq:cauchy} is the \emph{Cauchy
integral formula} up to factor ${2\pi} {\sqrt{ 1-\modulus{a}^2 }} $. Thus
we will write $f(a)={\left({2\pi \sqrt{ 1-\modulus{a}^2 }}\right)}^{-1}
-\widehat{f}(-a) $ for analytic extension of $f(\phi)$ to the unit disk.
The factor $2\pi$ is due to our normalization~\eqref{eq:n-measure} and
$\sqrt{ 1-\modulus{a}^2 }$ is connected with the invariant measure on $
\Space{D}{} $.

Consider now a realization of inverse formula~\eqref{eq:a-inverse}:
\begin{eqnarray}
f(e^{i\phi}) & = & \int_{ \Space{D}{} } \widehat{f}(a) \frac{\sqrt{1- 
\modulus{a}^2 }}{\bar{a}e^{i\phi}+1}\,d\mu(a) \nonumber \\
& = &-\int_{ \Space{D}{} } {2\pi }{\sqrt{ 1-\modulus{a}^2 }}{f}(-a) 
\frac{\sqrt{1- \modulus{a}^2 }}{\bar{a}e^{i\phi}+1}\,\frac{da}{{(1- 
\modulus{a}^2)}^2  } \nonumber \\
& = & {2\pi } \int_{ \Space{D}{} }
\frac{f(a)\,da}{(1-\bar{a}e^{i\phi}){(1-\modulus{a}^2)}^2}. 
\label{eq:ad-inverse}
\end{eqnarray}
The last integral is divergent, the singulyarity is concentrated near the 
unit circle. A regularization of the integral gives us the \emph{Sz\"ego} 
projection:
\begin{displaymath}
f(e^{i\phi}) = {\pi } \int_{\Gamma} \frac{f(a)\,da}{1-\bar{a}e^{i\phi}}.
\end{displaymath}
Values of $f(e^{i\phi})$ could be alternatively reconstructed from $f(a)$ 
by means of the limiting procedure $f(e^{i\phi})=\lim_{r \rightarrow 1 }
f(re^{i\phi})$.  Formula~\eqref{eq:ad-inverse} is of minor use in complex 
analysis but we will need it in Section~\ref{se:calculus} to construct
analytic functional calculus.
\end{example}  

\begin{example}
Let us consider second opportunity $X= \Space{D}{}$. For any integer 
$m\geq 2$ one could select a measure $d\mu_m(w)=4^{1-m}{(1- 
\modulus{w}^2 )}^{m-2} dw$, where $dw$ is the standard Lebesgue measure
on $ \Space{D}{} $, such that action~\eqref{eq:fr-lin} could be turn to 
a unitary representation~\cite[\S~IX.3]{Lang85}. Namely
\begin{equation}
[\pi_m(g) f](w)=f\left(  \frac{\alpha w + \beta}{\bar{\beta} w + 
\bar{\alpha}} \right) {(\bar{\beta} w + \bar{\alpha} )}^{-m}.
\end{equation}
If we again select $f_0(w)\equiv 1$ then
\begin{displaymath}
[\pi_m(g) f_0](w)= {(\bar{\beta} w + \bar{\alpha} )}^{-m},
\end{displaymath}
particularly $[\pi_m(f_\psi) f_0](w)= e^{im\psi} f_0(w)$ so this is 
again a vacuum vector with respect to $H$. The irreducible subspace $ 
\FSpace{F}{2}( \Space{D}{} )$ generated by $f_0(w)$ consists of analytic
functions and  is the $m$-th Bergman space (actually \person{Bergman}
considered only $m=2$). Now the transformation~\eqref{eq:berg-cauch} 
takes the form
\begin{eqnarray*}
\widehat{f}(a) & = & \scalar{f(w)}{[\pi_m(s(a)) f_0](w)}\\
& = & \left(1- \modulus{a}^2 \right)^{m/2} \int_{ \Space{D}{} }
\frac{f(w)}{{(a\bar{w}+1)}^m} \frac{dw}{{(1- \modulus{w}^2 )}^{2-m}},
\end{eqnarray*}
which is for $m=2$ the classical Bergman formula up to factor $\left(1-
\modulus{a}^2 \right)^{m/2}$.  Note that calculations in standard
approaches is ``rather lengthy and must be done in
stages''~\cite[\S~1.4]{Krantz82}.  As was mentioned early, its almost
identical with its inverse~\eqref{eq:a-inverse}:
\begin{eqnarray*}
f(w) & = & \int_{ \Space{D}{} } \widehat{f}(a) f_a(w)\,d\mu(a) = \int_{ 
\Space{D}{}
} \frac{\widehat{f}(a)}{{(1+\bar{a}w)}^m} \frac{da}{{(1- \modulus{a}^2
)}^2}
\\& = & \int_{ \Space{D}{} } f(a) \left( \frac{\sqrt{1- 
\modulus{a}^2 }}{1+\bar{a}w} \right)^m \frac{da}{{(1- \modulus{a}^2 
)}^2}
\\
& = & \int_{ \Space{D}{} }  \frac{f(a)}{(1+\bar{a}w) ^m} \frac{da}{{(1- 
\modulus{a}^2 )}^{2-m}},
\end{eqnarray*}
where $f(a)=\left(1-\modulus{a}^2\right)^{m/2} \widehat{f}(-a) $ is 
the same function as $f(w)$.
\end{example}  

\begin{example} \label{ex:segal}
The purpose of this Example is many-folds. First we would like to
justify the Segal-Bargmann name for the formula~\eqref{eq:a-inverse}.
Second we need a demonstration as our scheme will work for 
\begin{enumerate}
\item The nilpotent Lie group $ G=\Space{H}{n} $.
\item Non-geometric action of the group $G$.
\item Non-compact subgroup $H$.
\end{enumerate}

We will consider a representation of the Heisenberg group $ \Space{H}{n}
$ (see Section~\ref{se:manin}) on $ \FSpace{L}{2}( \Space{R}{n} ) $ by
operators of shift and multiplication~\cite[\S~1.1]{MTaylor86}:
\begin{equation} \label{eq:schrodinger}
g=(t,z): f(x) \rightarrow [\pi_{(t,z)}f](x)=e^{i(2t-\sqrt{2}qx+qp)} 
f(x- \sqrt{2}p), \qquad z=p+iq,
\end{equation}
i.e., this is the Schr\"odinger representation with parameter $\hbar=1$. 
As a subgroup $H$ we select the center of \Space{H}{n} consisting of 
elements $(t,0)$. It is non-compact but using the special form of 
representation~\eqref{eq:schrodinger} we could consider the 
cosets\footnote{$ \widetilde{G}$ is sometimes called the \emph{reduced} 
Heisenberg group. It seems that $ \widetilde{G}$ is a virtual object, 
which is important in connection with a selected representation of $G$.}
$ \widetilde{G} $ and $ \widetilde{H} $ of $G$ and $H$ by subgroup of 
elements $(\pi m,0)$, $m\in \Space{Z}{}$. Then~\eqref{eq:schrodinger}
also defines a representation of $ \widetilde{G} $ and $ \widetilde{H}
\sim \Gamma $. We consider the Haar measure on $ \widetilde{G} $ such
that its restriction on $ \widetilde{H} $ has the total mass equal to
$1$.

As ``vacuum vector'' we will select the original \emph{vacuum
vector} of quantum mechanics---the Gauss function $f_0(x)=e^{-x^2/2}$. 
Its transformations
are defined as follow:
\begin{eqnarray*}
f_g(x)=[\pi_{(t,z)} f_0](x) & = & e^{i(2t-\sqrt{2}qx+qp)}\,e^{-{(x- 
\sqrt{2}p)}^2/2}\\
     & = & e^{2it-(p^2+q^2)/2} e^{- ((p-iq)^2+x^2)/2+\sqrt{2}(p-iq)x} 
\\
     & = & e^{2it-z\bar{z}/2}e^{- (\bar{z}^2+x^2)/2+\sqrt{2}\bar{z}x}.
\end{eqnarray*}
Particularly $[\pi_{(t,0)} f_0](x)=e^{-2it}f_0(x)$, i.e., it really is 
a vacuum vector in the sense of our definition with respect to 
$\widetilde{H}$. Of course $\Omega=\widetilde{G} / \widetilde{H}$
isomorphic to  $\Space{C}{n}$ and mapping $s:
\Space{C}{n} \rightarrow \widetilde{G}$ simply is defined as $s(z)=(0,z)$. 
The Haar measure on $ \Space{H}{n} $ is coincide with the standard
Lebesgue measure on $ \Space{R}{2n+1} $~\cite[\S~1.1]{MTaylor86} thus
the invariant measure on $ \Omega $ also coincide with the Lebesgue
measure on $\Space{C}{n}$. Note also that composition law
$s^{-1}(g\cdot s(z))$ reduces to Euclidean shifts on $ \Space{C}{n} $. 
We also find $s^{-1}((s(z_1))^{-1}\cdot s(z_2))=z_2-z_1$ and 
$r((s(z_1))^{-1}\cdot s(z_2))= \frac{1}{2} \Im \bar{z}_1z_2$.

Transformation~\eqref{eq:berg-cauch} takes the form of embedding
$\FSpace{L}{2}(\Space{R}{n} ) \rightarrow \FSpace{L}{2}( \Space{C}{n} ) 
$ and is given by the formula
\begin{eqnarray}
\widehat{f}(z)&=&\scalar{f}{\pi_{s(z)}f_0}\nonumber \\
     &=&\pi^{-n/4}\int_\Space{R}{n} f(x)\, e^{-z\bar{z}/2}\,e^{-
(z^2+x^2)/2+\sqrt{2}zx}\,dx \nonumber \\
     &=&e^{-z\bar{z}/2}\pi^{-n/4}\int_\Space{R}{n} f(x)\,e^{-
(z^2+x^2)/2+\sqrt{2}zx}\,dx, \label{eq:tr-bargmann}
\end{eqnarray}
where $z=p+iq$. Then $\widehat{f}(g)$ belongs to
$\FSpace{L}{2}( \Space{C}{n} , dg)$ or its preferably to say that 
function $\breve{f}(z)=e^{z\bar{z}/2}\widehat{f}(t_0,z)$ belongs
to space $\FSpace{L}{2}( \Space{C}{n} , e^{- \modulus{z}^2 }dg)$ 
because $\breve{f}(z)$ is analytic in $z$. Such functions form the
\emph{Segal-Bargmann space}~\cite{Bargmann61,Segal60}  
$ \FSpace{F}{2}( \Space{C}{n}, e^{-
\modulus{z}^2 }dg) $ of functions, which are analytic by $z$ and
square-integrable with respect the Gaussian measure $e^{-
\modulus{z}^2}dz$. Analyticity of $\breve{f}(z)$ is equivalent to 
condition
$( \frac{ \partial }{  \partial\bar{z}_j } + \frac{1}{2} z_j I ) 
\widehat{f}(z)=0 $.

The integral in~\eqref{eq:tr-bargmann} is the well known 
Segal-Bargmann transform~\cite{Bargmann61,Segal60}. Inverse to it is 
given by a realization of~\eqref{eq:a-inverse}:
\begin{eqnarray}
f(x) & = & \int_{ \Space{C}{n} } \widehat{f}(z) f_{s(z)}(x)\,dz 
\nonumber\\
& = & \int_{ \Space{C}{n} } \breve{f}(z)  e^{- 
(\bar{z}^2+x^2)/2+\sqrt{2}\bar{z}x}\, e^{- \modulus{z}^2}\, dz
\end{eqnarray}
and this gives to~\eqref{eq:a-inverse} the name of Segal-Bargmann 
inverse. The corresponding operator $\oper{P}$~\eqref{eq:szego} is an 
identity operator $ \FSpace{L}{2}(\Space{R}{n}) \rightarrow 
\FSpace{L}{2}(\Space{R}{n}) $ and~\eqref{eq:szego} give an integral 
presentation of the Dirac delta. 

Meanwhile the orthoprojection  
$ \FSpace{L}{2}( \Space{C}{n},  e^{- \modulus{z}^2 }dg)  \rightarrow 
\FSpace{F}{2}( \Space{C}{n},  e^{- \modulus{z}^2 }dg) $ is of interest 
and is a principal ingredient in Berezin 
quantization~\cite{Berezin88,Coburn94a}. We could easy find its kernel 
from~\eqref{eq:s-b-proj}. Indeed, $ \widehat{f}_0(z)=e^{ - \modulus{z}^2 
} $, then kernel is 
\begin{eqnarray*}
K(z,w) & = & \widehat{f}_0(z^{-1}\cdot w) \bar{\chi}(r(z^{-1}\cdot w))\\
& = & \widehat{f}_0(w-z)e^{i\Im(\bar{z}w)} \\
& = & \exp(\frac{1}{2}(- \modulus{w-z}^2 +w\bar{z}-z\bar{w}))\\
& = &\exp(\frac{1}{2}(- \modulus{z}^2- \modulus{w}^2) +w\bar{z}).
\end{eqnarray*}
To receive the reproducing kernel for functions 
$\breve{f}(z)=e^{\modulus{z}^2} \widehat{f}(z) $ in the Segal-Bargmann 
space we should multiply $K(z,w)$ by $e^{(-\modulus{z}^2+ 
\modulus{w}^2)/2}$
which gives the standard reproducing kernel $= \exp(- \modulus{z}^2 
+w\bar{z})$ \cite[(1.10)]{Bargmann61}.
\end{example}

\begin{rem} \label{re:disk}
Started from the Hardy space $ \FSpace{H}{2}$ on $ \Gamma $ we 
recover the unit disk $ \Space{D}{} $ in a spirit of the Erlangen
program as a homogeneous space connected with its group of symmetry. On
the other hand $ \Space{D}{} $ carries out a good portion  of
information about space of maximal ideal of $ \FSpace{H}{\infty} $,
which is a subspace of $ \FSpace{H}{2}$ forming an algebra under
multiplication~\cite[Chap.~6]{Douglas72a}. So we meet $ \Space{D}{} $ 
if we will look for geometry
of $ \FSpace{H}{p} $ from the Descartes viewpoint. This illustrate once
more that there is an intimate connection between these two approaches.
\end{rem}
\begin{rem}
The abstract scheme presented here is not a theory with the only
example. One could apply it to analytic function theory of several
variables. This gives a variety of different function theories both known 
(see below) and new~\cite{CnopsKisil97a,Kisil97c}. For
example, the biholomorphic automorphisms~\cite[Chap.~2]{Rudin80} of unit
ball in $ \Space{C}{n} $ lead to several complex variable theory, while
M\"obius transformations~\cite{Cnops94a} of $ \Space{R}{n} $ guide to
Clifford analysis~\cite{BraDelSom82}. It is interesting
that both types of transformations could be represented as
fraction-linear ones associated with some $2 \times 2 $-matrixes and 
very reminiscent the content of this Section. 

On the other hand, historically several complex variable theory was
developed as coordinate-wise extension of notion of holomorphy while
Clifford analysis took a conceptually different route of a factorization
of the Laplacian. It is remarkable that such differently rooted
theories could be unified within presented scheme. Moreover, relations 
between symmetry groups of complex and Clifford analysis explain why one
could make conclusions about complex analysis using Clifford technique. 
This repeats relations between affine and Euclidean geometry, for 
example.
\end{rem}  
Considering such a pure mathematic theory as complex analysis we 
continuously use the language of applications (vacuum vector, coherent 
states, wavelet transform, etc.). We will consider these applications 
explicitly later in Section~\ref{se:observ}.

\section{Non-Commutative Conformal Geometry} \label{se:calculus}
\epigraph{The search for generality and unification is one of the distinctive 
features of twentieth century mathematics, and functional analysis seeks to 
achieve these goals.}{M.~Kline \cite[\S~46.1]{Kline72}.}{}

Accordingly to \person{F.~Klein} one could associate to every
specific geometry a group of symmetries, which
characterizes it.  This is also true for new geometries discovered in 20-th 
century. Many books on functional analysis (see for 
example~\cite{Douglas72a,KirGvi82}) contain chapters entitled as
\emph{Geometry of Hilbert Spaces} where they present inner product,
Pythagorean theorem, orthonormal basis etc. All mentioned notions are
invariants of the \emph{group of unitary transformation} of a Hilbert
space. On the contrary, the geometry of Hilbert spaces is not directly
related to any coordinate algebra. So this is a geometry more in Klein
meaning than in Descartes ones.

Considering Banach spaces (or more generally locally convex ones) one
has group of continuous \emph{affine} transformations. Its invariant are
points, lines, hyperplanes and notions, which could be formulated in
their term. The excellent example is \emph{covexity}. Being continuous
these transformations have topological invariants like compactness.
Thereafter the Krein-Milman theorem is formulated entirely in terms of
continuous affine invariants and thus belong to the body of
\emph{infinite dimensional continuous affine geometry}.

Which groups of symmetries produce an interesting geometries in Banach
algebras? One could observe that $\SL $ defines a
fraction-linear transformation
\begin{equation} \label{eq:t-fr-lin}
g: t \mapsto    
g \cdot t={(\bar{ \beta}t +\bar{\alpha})^{-1}}{(\alpha t + \beta)},
\qquad g^{-1}=\matr{\alpha}{\beta}{\bar{\beta}}{\bar{\alpha}}.           
\end{equation}
not only for $t$ being a complex number but also any element of a Banach
algebra $ \algebra{A} $ provided $(\bar{ \beta}t +\bar{\alpha})^{-1}$
exists. To be sure that \eqref{eq:t-fr-lin} is defined for any $g\in
\SL$ one could impose the condition $ \norm{t} < 1$. Let
us fix such a $t\in \algebra{A}$, $ \norm{t}\leq 1$. Then elements
$g \cdot t$, $g\in \SL$ form a subset $ \Space{T}{}
$ of $ \algebra{A} $ invariant under $\SL$ and is a
homogenous space where $\SL$ acts by the rule:
$g: {g'}\cdot t \mapsto {(gg')}\cdot t$. Note also that 
$ \norm{ g\cdot t} < 1$. Thus it is reasonably to consider $
\Space{T}{} $ as an analog of the unit disk $ \Space{D}{} $ and try to
construct analytic functions on $ \Space{T}{} $. This problem is known
as a functional calculus of an operator and we would like to present its
solution inspired by the Erlangen program.

We are looking for a possibility to assign an element $\Phi(f,{s(a)} 
\cdot t)=f({s(a)}\cdot t)$ in $ \algebra{A} $ to every pair $f\in 
\FSpace{F}{2}( \Gamma ) $ and $ {s(a)}\cdot t \in \Space{T}{} $, $a\in 
\Space{D}{} $.
If we consider $\FSpace{F}{2}( \Gamma ) $ as a linear space then 
naturally to ask that mapping $\Phi: (f,{s(a)}\cdot t) \rightarrow 
f({s(a)}\cdot t)$ should be linear in $f$.
In line with function theory~\eqref{eq:g-transform} one also could define 
associated representation of $ \SL $ by the rule:
\begin{equation} \label{eq:t-transform}
\tau_{g} \Phi (f,g'\cdot t)=(\bar{\beta}t+\bar{\alpha})^{-1}
\Phi(f,{(gg')}\cdot t).
\end{equation}
And finally one could assume that~\eqref{eq:g-transform} 
and~\eqref{eq:t-transform} are in the agreement:
\begin{equation} \label{eq:intertwin}
\tau_{g} \Phi (f,g'\cdot t)=\Phi (\pi_{g} f ,g'\cdot t).
\end{equation}
Thereafter one could obtain the functional calculus $\Phi(f, g \cdot t)$
from a knowledge of $\Phi(f,  t)$ for all $f \in \FSpace{F}{2}( \Gamma 
)$ using~\eqref{eq:t-transform} and~\eqref{eq:intertwin}:
\begin{displaymath}
\Phi(f, g \cdot t)= (\bar{\beta}t+\bar{\alpha}) \Phi(\pi_g f, t).
\end{displaymath}
However it is easier to process in another way around. One could define 
$\Phi(f_0, g \cdot t)$ for a specific $f_0$ and all $g$ and then extend 
it to an arbitrary $f$ using~\eqref{eq:intertwin} and linearity. We 
are now ready to give 
\begin{defn}
The \emph{functional calculus} $\Phi(f, g\cdot t)$ for an element $t$ of 
a Banach algebra $ \algebra{A} $ is a linear continuous map $\Phi: 
\FSpace{H}{2} \times \Space{T}{} \rightarrow \algebra{A} $ satisfying to
the following conditions 
\begin{enumerate} 
\item $\Phi$ intertwines two representations $\pi_g$ 
\eqref{eq:g-transform} and $\tau_g$ \eqref{eq:t-transform} of $\SL$,
namely~\eqref{eq:intertwin} holds. \item $\Phi(f_0,g\cdot t)=e$ for all
$g$, where $f_0(\phi)\equiv 1$ and $e$ is the unit of $ \algebra{A} $.
\end{enumerate} 
\end{defn}
We know~\eqref{eq:ad-inverse} that any function $f(e^{i\phi})$ has a 
decomposition
\begin{displaymath}
f(e^{i\phi})= 2\pi\int_{ \Space{D}{} } 
\frac{f(a)\,da}{(1-\bar{a}e^{i\phi})(1-\modulus{a}^2)}.
\end{displaymath}
Thus using linearity
\begin{eqnarray}
\Phi(f(e^{i\phi}),t) & = & \Phi(2 \pi \int_{ \Space{D}{} } 
\frac{f(a)\,da}{(1-\bar{a}e^{i\phi})(1-\modulus{a}^2)}\,,\, t) \nonumber 
\\
 & = & 2\pi \int_{ \Space{D}{} } f(a)
\Phi(\frac{1}{1-\bar{a}e^{i\phi}}, t) \frac{da}{ (1-\modulus{a}^2)} 
\label{eq:use-lin} \\ 
& = & -2\pi \int_{ \Space{D}{} } f(a) \Phi(\pi_{s(-a)} f_0, t)
\frac{da}{ 1-\modulus{a}^2} \nonumber \\
& = & -2\pi \int_{ \Space{D}{} } f(a) \tau_{s(-a)}\Phi( f_0, t)
\frac{da}{ 1-\modulus{a}^2} \label{eq:use-inter} \\
& = & 2 \pi \int_{ \Space{D}{} } f(a) (1-\bar{a}t)^{-1}\Phi( f_0, 
s(-a)\cdot t) \frac{da}{ 1-\modulus{a}^2} \nonumber \\
& = & 2 \pi \int_{ \Space{D}{} } f(a) (1-\bar{a}t)^{-1}  \frac{da}{ 
1-\modulus{a}^2} \label{eq:t-bergman}.
\end{eqnarray}
Here we use linearity of functional calculus to
receive~\eqref{eq:use-lin}, its intertwining property to
get~\eqref{eq:use-inter}, and finally  \eqref{eq:t-bergman} that $e$ is 
always the image of $1$. The only difference between the last
formula~\eqref{eq:t-bergman}
and~\eqref{eq:ad-inverse} 
is the placing of $t$ instead of $e^{i\phi}$. If one would like to 
obtain an integral formula for the functional calculus in terms of 
$f(e^{i\phi})$ itself rather than $f(a)$ then
\begin{eqnarray}
\Phi(f(e^{i\phi}),t) & = & {\pi } \int_{ \Space{D}{} } f(a) 
(1-\bar{a}t)^{-1}  \frac{da}{ 1-\modulus{a}^2} \nonumber \\
& = & {\pi } \int_{ \Space{D}{} } \int _{ \Gamma }   
\frac{f(e^{i\phi})\, de^{i\phi}}{a-e^{i\phi}} (1-\bar{a}t)^{-1} \frac{da}{
1-\modulus{a}^2} \nonumber \\ 
& = & {\pi }  \int _{ \Gamma }   
f(e^{i\phi}) \int_{ \Space{D}{} } \frac{(1-\bar{a}t)^{-1}}{a-e^{i\phi}}  
\frac{da}{ 1-\modulus{a}^2}\, de^{i\phi}\nonumber \\ 
& = & \int _{ \Gamma } f(e^{i\phi}) (t-e^{i\phi})^{-1} 
de^{i\phi}
\end{eqnarray}
From our explicit construction it follows that functional calculus is 
\emph{always exist} and \emph{unique} under assumption $ \norm{t}<1 $. 
Moreover from the last formula it follows also that our functional
calculus coincides with the Riesz-Dunford
one~\cite[\S~VIII.2]{DunfordSchwartzI}.
\begin{rem} \label{re:calc}
Guided by the Erlangen program we constructed a functional calculus
coincided with the classical Dunford-Riesz calculus. The former is
traditionally defined in terms of algebra homomorphisms, i.e., on the
Descartes geometrical language. So two different approaches again lead
to the same answer.
\end{rem} 
If one compares epigraphs to Sections~\ref{se:analytic}
and~\ref{se:calculus} then it will be evidently that the functional 
calculus based on group representations responds to ``the search for 
generality and unification''.
Besides a funny alternative to classic Riesz-Dunford calculus the 
given scheme give unified approach to a variety of existing and yet to 
be discovered calculi~\cite[Remark~4.2]{Kisil95i}. We consider
some of them in next Section in connections with physical applications.

\section{Do We Need That Observables Form an Algebra?}\label{se:observ}
\epigraph{\ldots{} Descartes himself thought that all of
physics could be reduced to geometry.}{M.~Kline
\cite[\S~15.6]{Kline72}.}{}

We have already mentioned that the development of non-commutative 
geometry intimately connected with physics, especially with quantum 
mechanics. Despite the fact that quantum mechanics attracted
considerable efforts of physicists and mathematicians during this 
century we still have not got at hands a consistent quantum theory.
Unsolved paradoxes and discrepant interpretations are linked with very
elementary quantum models. As a solution to this frustrating situation
some physicists adopt a motto that science is permitted to supply
inconsistent models ``as is'' so far they give reasonable numerical
predictions to particular experiments\footnote{And even more strange:
there exists a variety of ``immortal'' theories already contradicting to
our experience, see for example the red shift problem
in~\cite{Segal90a}.}.
However other researchers still believe that there exists a difference 
between the standards of scientific enterprise and software selling 
industry.  

We will review some topics in quantization problem in connection with 
the two approaches to non-commutative geometry. Quantum mechanics 
begun from experiments clearly indicated that Minkowski four
dimensional space-time could not describe microscopical structure of our
world (on the other hand it is also unsuitable for macroscopic 
description and these two faults became connected~\cite{Segal90a}). A 
way out is searched in many directions: dimensionality of
microscopical space was raised to $m>4$ dimensions, its 
discontinuity was imposed, etc. These search is based on the following 
conclusion~\cite{Segal96b}:
\begin{quote}
From a quantum mechanical standpoint, the trust of these works was that 
``space-time'' was more logically not a primary concept but one derived
from the algebra of ``observables'' (or operators). In its simplest form
the idea was space-time might be, at a deeper level, the spectrum of an 
appropriate commutative subalgebra that is invariant under the 
fundamental symmetry group.
\end{quote}
The right algebra of quantum observables is often constructed from 
classic observables by means of \emph{quantization}.
A summary of quantization problem could be given as follows (see 
also~\cite[\S~15.4]{Kirillov76} and Table~\ref{ta:c-q}). The
principal objects of classical mechanics in Hamiltonian formalism are:
\begin{enumerate}
\item The \emph{phase space}, which is a smooth symplectic manifold $M$.
\item \emph{Observables} are real functions on $M$ and could be 
considered as elements of a linear space or algebra.
\item \emph{States} of a physical system are linear functionals on 
observables.
\item \emph{Dynamics} of observables is defined by a Hamiltonian $H$ and
equation
\begin{equation} \label{eq:hamilton}
\dot{F}=\{H,F\}.
\end{equation}
\item \emph{Symmetries} of the physical system act on observables or 
states via canonical transformations of $M$.
\end{enumerate}
Of course points of phase space is in one-to-one correspondence with 
states of the system and they together are completely defined by the set 
of observables. The later viewed as an algebra recovers $M$ as its space
of maximal ideals. Thus observables, dynamics, and symmetries are 
primary objects while phase space and states could be restored from
them.

The quantum mechanics has its parallel set of notions:
\begin{enumerate}
\item \emph{Phase space} is the projective space $P(V)$ of a Hilbert 
space $V$.
\item \emph{Observables} are self-adjoint operators on $V$.
\item \emph{States} of a physical system defined by a vector 
$\xi \in V$, $ \norm{\xi}=1 $.
\item \emph{Dynamics} of an observable $ \widehat{F}$ defined by a 
self-adjoint operator $ \widehat{H}$ via the Heisenberg equation
\begin{equation}
\dot{ \widehat{F} }= \frac{i \hbar }{2\pi} [ \widehat{H}, \widehat{F} ].
\end{equation}
\item \emph{Symmetries} of physical system act on observables or 
states via unitary operators on $V$.
\end{enumerate}
Again we could start from an algebra of observables and then realize states
and Hilbert space $V$ via GNS construction~\cite[\S~4.3]{Kirillov76}. 

The similarity of two descriptions (with correspondence of the Poisson 
bracket and commutator) makes it very tempting to construct
a map $\oper{Q}$, which is called the \emph{Dirac quantization}, such 
that
\begin{enumerate}
\item $\oper{Q}$ maps a function $F$ to an operator $ 
\widehat{F}= \oper{Q}(F) $, i.e., classic observables to quantum ones.
\item $\oper{Q}$ has an algebraic property $ \widehat{\{F_1,F_2\}}=[ 
\widehat{F}_1, \widehat{F}_2 ] $, i.e., maps the Poisson brackets to the 
commutators.
\item $\oper{Q}(1)=I$, i.e., the image of the function $1$ identically 
equal to one is the identity operator $I$.
\end{enumerate}
Unfortunately, $\oper{Q}$ does not exists even for very simple cases. 
For example, the ``no go'' theorem of \person{van Hove} 
in~\cite{Folland89} states that the canonical commutation 
relations~\eqref{eq:a-heisenberg} could not be extended even for 
polynomials in $X$ and $Y$ of degree greater than $2$.

As a solution one could try to select a primary set of classic
observables $F_1$, \ldots, $F_n$, which defines coordinates on $M$ and 
thus any other observable $F$ is a function of $F_1$, \ldots, $F_n$: 
$F=F(F_1,\ldots,F_n)$. For such a small set of observables a 
quantization $\oper{Q}$ always exists (for example the \emph{geometric 
quantization}~\cite[\S~15.4]{Kirillov76}) and one could hope to 
construct general quantum observable $ \widehat{F} $ as a function of 
the primary ones $ \widehat{F}=F( \widehat{F}_1, \ldots, \widehat{F}_n ) 
$.

But to do that one should answer the question equivalent to functional 
calculus: \emph{What is a function of $n$ operators}? For a commuting
set of operators as well as for the case $n=1$ it is easy to give an
answer in terms of \emph{algebra homomorphisms}: the functional calculus
$\Phi$ is an algebra homomorphism from an algebra of function to an
algebra of operators, which is fixed by a condition $\Phi(F_i)=
\widehat{F}_i $, $i=1,\ldots,n$. But for physically interesting case of
several non-commuting $ \widehat{F}_i$ this way is impossible. 

Several solutions were given whether analytically by integral formulas
like for the Feynman~\cite{Feynman51} and Weyl~\cite{Anderson69} calculi
or algebraically by ``ordering rules'' like for the
Wick~\cite[\S~6]{Segal96a} and anti-Wick~\cite{Berezin88} calculi. These
approaches are connected: if calculus is defined via integral formula
then it correspond to some ordering (e.g., the Feynman case) or
symmetrization (e.g., the Weyl case); and \emph{vise versa} starting 
from particular ordering a suitable integral representation could be
found~\cite{Berezin88}. 

We will consider the Weyl functional calculus in more details. 
Let $\mathbf{T}=\{ T_j\} $, $j=1,\ldots,n$ be a set of 
bounded selfadjoint operators. Then for every set $\{ \xi_j\}$ of real 
numbers the unitary operator $\exp(i\sum_1^n \xi_j T_j)$ is well 
defined.
The \emph{Weyl functional calculus}~\cite{Anderson69} is defined by the 
formula
\begin{equation}
f(\mathbf{T})=\int_{ \Space{R}{n} } \exp \left(i\sum_1^n \xi_j T_j 
\right) \widehat{f}(\xi)\, d\xi,
\end{equation}
where $\widehat{f}(\xi)$ is usual Fourier transform of $f(x)$. 
In~\cite{Anderson69} the following properties of the Weyl calculus were 
particularly proved:
\begin{enumerate}
\item The image of a polynomial $p(x)$ is the symmetric in $T_j$ 
polynomial $p(\mathbf{T})$ (e.g. image of $p(x_1,x_2)=x_1x_2$ is 
$p(T_1,T_2)= \frac{1}{2}(T_1 T_2 + T_2 T_1) $).
\item The Weyl calculus commutes with affine transformations, i.e., let 
$S_i=\sum_{j=1}^n m_{ij} T_j$ for a nonsingular matrix $\{ m_{ij}\}$, 
then 
\begin{equation}  \label{eq:w-inter}
f(\mathbf{S})=g(\mathbf{T}) \textrm{ where  } g(x)=f(m\cdot x).
\end{equation}
\end{enumerate}
Moreover the last condition together with the requirement that Weyl 
functional calculus gives the standard functional calculus for 
polynomials in one variable completely defines it. We know that 
functional calculus for one variable could be also defined by its 
covariance property. Thus we could give an equivalent
definition\footnote{We skip details because the problem is beyond the $
\FSpace{L}{2} $ scope presented in the paper.}
\begin{defn}
The Weyl functional calculus for an $n$-tuple $\mathbf{T}$ of selfadjoint
operators is a linear map from a space of functions on
$ \Space{R}{n}$ to operators in a \Cstar-algebra $ \algebra{A} $ such 
that
\begin{enumerate}
\item Intertwines two representation of affine transformation, namely 
\eqref{eq:w-inter} holds.
\item Maps the function $e^{ix_1}$ on $ \Space{R}{} $ to operator 
\begin{equation} \label{eq:exponent}
e^{iT_1} = \sum \frac{(iT_1)^n}{n!}.
\end{equation}
\end{enumerate}
\end{defn}  

While the affine group is natural for Euclidean space one could need 
other groups under different environments. For example, one could
consider a non-commutative geometry generated by
the group of conformal (or M\"obius) transformations. It is
known~\cite{Cnops94a} that the group $M$ of M\"obius (mapping spheres 
onto spheres) transformations of $\Space{R}{n}$ could be represented via
fraction-linear transformations with the help of Clifford
algebras $\Cliff{n}$.  Taking the tensor product $ 
\widetilde{\algebra{A}}=\algebra{A} \otimes \Cliff{n}$ of an operator
algebra $\algebra{A}$ and Clifford algebra we could
consider a representation of conformal group $M$ in the
non-commutative space $\widetilde{\algebra{A}}$, which also is defined 
by fraction-linear transformations.

Such an approach leads particularly to a natural definition
of joint spectrum of $n$-tuples of non-commuting operators
and a functional calculus of them~\cite{Kisil95i}.  The functional
calculus is not an algebra homomorphism (for $n>1$) but is an
intertwining operator between two representations of
conformal group $M$. Nevertheless the joint spectrum obeys
a version of the spectral mapping theorem~\cite{Kisil95i}.

It is naturally to expect that physically interesting functional calculi
from group representations are not restricted to two given examples.
Looking on these functional calculi as on a quantization one
found that common contradiction between the group covariance and an
algebra homomorphism property vanishes (as the algebra homomorphism
property itself). It turns that the group covariance alone is sufficient
for the definition of functional calculus. Moreover such an approach 
gives a unified description of all known functional calculi and suggests 
new ones.

One should be aware that algebraic structure it still important: it
enters to the definition of the group action via fraction-linear
transformations~\eqref{eq:t-transform} or definition of 
exponent~\eqref{eq:exponent}. So
the question posed at the title of Section might
have an answer: \emph{Yes, we need that observables form an algebra, 
but we do not need that quantization will be connected with an algebra 
homomorphism property in general.}

\section{Conclusion}
\epigraph{We often hear that mathematics consists mainly in ``proving 
theorems''. Is a writer's job mainly that of ``writing
sentences''?}{G.-C.~Rota \cite[Chap.~11]{KacRotaSchwartz92}.}{}
Looking through the present discussion one could not resist to the 
following feeling:

\emph{The goal of new mathematical facts is to elaborate the true 
language for an understanding old ideas, such that they will be 
expressed with the ultimate perfectness.}

This goal could not be obviously accomplished within a finite time. So 
we have a good reason for self-irony mentioned at the end of 
Introduction.

\section{Acknowledgments}
The paper was written while the author stay at the Department of
Mathematical Analysis, University of Gent, as a Visiting Postdoctoral
Fellow of the FWO-Vlaanderen (Fund of Scientific Research-Flanders),
Scientific Research Network ``Fundamental Methods and Technique in
Mathematics'' 3GP03196, whose hospitality and support I gratefully
acknowledge.

I would like to thank to Prof.~\person{I.E.~Segal} and \person{Yu.~Manin}
for useful remarks on the subject.  I am also grateful to
Dr.~\person{A.A.~Kirillov, Jr.} for critical discussion, which allows to
improve presentation in Section~\ref{se:manin}.

In the paper we try to establish some connections between several 
areas of algebra, analysis, and geometry. A reasonable bibliography to
such a paper should be necessary incomplete and the given one is
probably even not representative. Any suggestions about unpardonable
omitted references will be thankfully accepted.

\small
\bibliographystyle{plain}
\bibliography{MRABBREV,analyse,acombin,aphilos,aphysics}

\end{document}

%% file: mydef.tex
\def\ifundefined#1{\expandafter\ifx\csname#1\endcsname\relax}
\newcommand{\comment}[1]{}

\usepackage{amsfonts}

\newcommand{\algebra}[1]{\ensuremath{\mathfrak{#1}}}
\newcommand{\Heisen}[1]{\ensuremath{\mathbb{H}^{#1}}}

\newcommand{\Cliff}[2][\comment]{\ensuremath{{\cal C}\kern-0.18em\ell(#1,#2)}}

\newcommand{\Space}[2]{\ensuremath{ {\mathbb{#1}^{#2}} }}
\newcommand{\such}{\,\mid\,}

\newcommand{\FSpace}[2]{{\ensuremath{ #1_{#2} }}}

\ifundefined{qed}
    \DeclareMathSymbol{\qed}{0}{AMSa}{"03}
\fi
\newcommand{\Cstar}{$C^{*}$}

\newcommand{\norm}[1]{\left\| #1 \right\|}
\newcommand{\modulus}[1]{\left| #1 \right|}
\newcommand{\scalar}[2]{\langle #1,#2\rangle}

\providecommand{\eqref}[1]{\textup{(\ref{#1})}}

\newcommand{\person}[1]{\textsc{#1}}

%% file: nccg1.bbl
\newcommand{\noopsort}[1]{} \newcommand{\printfirst}[2]{#1}
  \newcommand{\singleletter}[1]{#1} \newcommand{\switchargs}[2]{#2#1}
  \newcommand{\irm}{\textup{I}} \newcommand{\iirm}{\textup{II}}
  \newcommand{\vrm}{\textup{V}} \providecommand{\noopsort}[1]{}
  \providecommand{\printfirst}[2]{#1} \providecommand{\singleletter}[1]{#1}
  \providecommand{\switchargs}[2]{#2#1}
\begin{thebibliography}{10}

\bibitem{AliAntGazMue}
S.~Twareque Ali, J.-P. Antoine, J.-P. Gazeau, and U.A. Mueller.
\newblock Coherent states and their generalizations: {A} mathematical overview.
\newblock {\em Rev. Math. Phys.}, 7(7):1013--1104, 1995.
\newblock Zbl \# 837.43014.

\bibitem{Anderson69}
Robert~F.~V. Anderson.
\newblock The {Weyl} functional calculus.
\newblock {\em J. Funct. Anal.}, 4:240--267, 1969.

\bibitem{AtiyahSchmid80}
Michael Atiyah and Wilfried Schmid.
\newblock A geometric construction of the discrete series for semisimple {Lie}
  group.
\newblock In J.A. Wolf, M.~Cahen, and M.~De Wilde, editors, {\em Harmonic
  Analysis and Representations of Semisimple {Lie} Group}, volume~5 of {\em
  Mathematical Physics and Applied Mathematics}, pages 317--383. D. Reidel
  Publishing Company, Dordrecht, Holland, {\noopsort{}}1980.

\bibitem{Bargmann61}
V.~Bargmann.
\newblock On a {H}ilbert space of analytic functions and an associated integral
  transform. {Part I}.
\newblock {\em Comm. Pure Appl. Math.}, 3:215--228, 1961.

\bibitem{Berezin88}
Felix~A. Berezin.
\newblock {\em Method of Second Quantization}.
\newblock ``Nauka'', Moscow, {\noopsort{}}1988.
\newblock (Russian).

\bibitem{BraDelSom82}
F.~Brackx, R.~Delanghe, and F.~Sommen.
\newblock {\em Clifford Analysis}, volume~76 of {\em Research Notes in
  Mathematics}.
\newblock Pitman Advanced Publishing Program, Boston, 1982.

\bibitem{ChamConnes96b}
Ali Chamseddine and Alain Connes.
\newblock The spectral action principle.
\newblock {\em Commun. Math. Phys}, 186(3):731--750, 1997.
\newblock preprint \texttt{hep-th/9606001}.

\bibitem{ChanRotaStein95}
Wendy Chan, Gian-Carlo Rota, and Joel~A. Stein.
\newblock The power of positive thinking.
\newblock In {\em Invariant Methods in Discrete and Computational Geometry.
  Proceedings of the conference held in {Cura\c{c}ao, June} 13--17, 1994},
  pages 1--36. Kluwer Academic Publishers, Dordrecht, {\noopsort{}}1995.

\bibitem{Cnops94a}
Jan Cnops.
\newblock {\em {Hurwitz} Pairs and Applications of {M\"obius} Transformations}.
\newblock {Habilitation} dissertation, Universiteit Gent, Faculteit van de
  Wetenschappen, 1994.
\newblock \texttt{ftp://cage.rug.ac.be/pub/clifford/jc9401.tex}.

\bibitem{CnopsKisil97a}
Jan Cnops and Vladimir~V. Kisil.
\newblock Monogenic functions and representations of nilpotent lie groups in
  quantum mechanics.
\newblock 1997.
\newblock (In preparation).

\bibitem{Coburn94a}
Lewis~A. Coburn.
\newblock {Berezin-Toeplitz} quantization.
\newblock In {\em Algebraic Mettods in Operator Theory}, pages 101--108.
  Birkh\"auser Verlag, New York, 1994.

\bibitem{Connes94}
Alain Connes.
\newblock {\em Non-Commutative Geometry}.
\newblock Academic Press, New York, {\noopsort{}}1994.

\bibitem{Daubechies92}
Ingrid Daubechies.
\newblock {\em Ten Lectures on Wavelets}, volume~61 of {\em CBMS-NSF Regional
  Conference Series in Applied Mathematics}.
\newblock Society for Industrial and Applied Mathematics (SIAM), Philadelphia,
  PA, {\noopsort{}}1992.

\bibitem{Douglas72a}
Ronald~G. Douglas.
\newblock {\em {B}anach Algebra Techniques in Operator Theory}, volume~49 of
  {\em Pure and Applied Mathematics}.
\newblock Academic Press, Inc., New York, {\noopsort{}}1972.

\bibitem{DunfordSchwartzI}
Nelson Dunford and Jacob~T. Schwartz.
\newblock {\em Linears Operators. Part I: General Theory}, volume VII of {\em
  Pure and Applied Mathematics}.
\newblock John Wiley \& Sons, Inc., New York, 1957.

\bibitem{Feynman51}
Richard~P. Feynman.
\newblock An operator calculus having applications in quantum electrodynamics.
\newblock {\em Phys. Rev. A (3)}, 84(2):108--128, 1951.

\bibitem{Folland89}
Gerald~B. Folland.
\newblock {\em Harmonic Analysis in Phase Space}.
\newblock Princeton University Press, Princeton, New Jersey, {\noopsort{}}1989.

\bibitem{GilbertMurray91}
John~E. Gilbert and Margaret~A.M. Murray.
\newblock {\em Clifford Algebras and {Dirac} Operators in Harmonic Analysis},
  volume~26 of {\em Cambridge Studies in Advanced Mathematics}.
\newblock Cambridge University Press, Cambridge, {\noopsort{}}1991.

\bibitem{Howe80a}
Roger Howe.
\newblock On the role of the {Heisenberg} group in harmonic analysis.
\newblock {\em Bull. Amer. Math. Soc. (N.S.)}, 3(2):821--843, 1980.

\bibitem{HoweTan92}
Roger Howe and Eng~Chye Tan.
\newblock {\em Non-Abelian Harmonic Analysis: Applications of
  ${SL(2,\Space{R}{})}$}.
\newblock Universitext. Springer-Verlag, New York, 1992.

\bibitem{KacRotaSchwartz92}
Mark Kac, Gian-Carlo Rota, and Jacob~T. Schwartz.
\newblock {\em Discrete Thoughts: Essays on Mathematics, Science, and
  Philosopy}.
\newblock Birkh\"auser Verlag, Boston, {\noopsort{1984}}1992.

\bibitem{Kassel94}
Christian Kassel.
\newblock {\em Quantum Groups}, volume 155 of {\em Graduate Text in
  Mathematics}.
\newblock Springer-Verlag, New York, {\noopsort{}}1994.

\bibitem{KirGvi82}
A.~A. Kirillov and A.~D. Gvishiani.
\newblock {\em Theorems and Problems in Functional Analysis}.
\newblock Problem Books in Mathematics. Springer-Verlag, New York,
  {\noopsort{}1982}.

\bibitem{Kirillov76}
Alexander~A. Kirillov.
\newblock {\em Elements of the Theory of Representations}, volume 220 of {\em A
  Series of Comprehensive Studies in Mathematics}.
\newblock Springer-Verlag, New York, 1976.

\bibitem{Kisil95a}
Vladimir~V. Kisil.
\newblock Integral representation and coherent states.
\newblock {\em Bull. Belg. Math. Soc. Simon Stevin}, 2(5):529--540, 1995.

\bibitem{Kisil95i}
Vladimir~V. Kisil.
\newblock {M\"obius} transformations and monogenic functional calculus.
\newblock {\em ERA of Amer. Math. Soc.}, 2(1):26--33, August 1996.

\bibitem{Kisil97c}
Vladimir~V. Kisil.
\newblock Analysis in{ $\Space{R}{1,1}$} or the principal function theory.
\newblock page~25, 1997.
\newblock E-print \texttt{funct-an/9712003}, \texttt{clf-alg/kisi9703}.

\bibitem{KlaSkag85}
John~R. Klauder and Bo-Sture Skagerstam, editors.
\newblock {\em Coherent States. Applications in physics and mathematical
  physics.}
\newblock World Scientific Publishing Co., Singapur, {\noopsort{}}1985.

\bibitem{Kline72}
Morris Kline.
\newblock {\em Mathematical Thought from Ancient to Modern Times}.
\newblock Oxford University Press, New York, 1972.

\bibitem{KnappWallach76}
A.W. Knapp and N.R. Wallach.
\newblock Szeg\"o kernels associated with discrete series.
\newblock {\em Invent. Math.}, 34(2):163--200, 1976.

\bibitem{TeXbook}
Donald~E. Knuth.
\newblock {\em The {\TeX book}}.
\newblock Addison-Wesley, Reading, Massachusetts, {\noopsort{1984}}1986.

\bibitem{Koranyi72a}
Adam Koranyi.
\newblock Harmonic functions on symmetric spaces.
\newblock In William~M. Boothby and Guido~L. Weiss, editors, {\em Symmetric
  Spaces: Short Courses Presented at {Washington} University}, volume~8 of {\em
  Pure and Applied Mathematics}, pages 379--412. Marcel Dekker, Inc, New York,
  {\noopsort{}}1972.

\bibitem{Krantz82}
Steven~G. Krantz.
\newblock {\em Function Theory of Several Complex Variables}.
\newblock Pure and Applied Mathematics. John Wiley \& Sons, New York,
  {\noopsort{}}1982.

\bibitem{Lang85}
Serge Lang.
\newblock {\em {$SL_2(\mathbf{R})$}}, volume 105 of {\em Graduate Text in
  Mathematics}.
\newblock Springer-Verlag, New York, {\noopsort{}}1985.

\bibitem{Manin88a}
Yuri~I. Manin.
\newblock {\em Quantum Group and Non-Commutative Geometry}.
\newblock Universit\'e de Monr\'eal, Centre de Recherches Math\'ematiques,
  Monr\'eal, PQ, {\noopsort{}}1988.

\bibitem{Manin91}
Yuri~I. Manin.
\newblock {\em Topics in Noncommutative Geometry}.
\newblock Princeton University Press, Princeton, New Jersey, {\noopsort{}}1991.

\bibitem{Perelomov86}
A.~M. Perelomov.
\newblock {\em Generalized Coherent States and Their Applications}.
\newblock Springer-Verlag, Berlin, {\noopsort{}}1986.

\bibitem{Rudin80}
Walter Rudin.
\newblock {\em Function Theory in the Unit Ball of{ \Space{C}{n}}}.
\newblock Springer-Verlag, New York, {\noopsort{}}1980.

\bibitem{Sabinin72}
Lew~I. Sabinin.
\newblock Loop geometries.
\newblock {\em Math. Notes}, 12:799--805, 1972.

\bibitem{Segal60}
Irving~E. Segal.
\newblock {\em Mathematical Problems of Relativistic Physics}, volume~II of
  {\em Proceedings of the Summer Seminar (Boulder, Colorado, 1960)}.
\newblock American Mathematical Society, Providence, R.I., 1963.

\bibitem{Segal90a}
Irving~E. Segal.
\newblock The mathematical implications of fundamental physical principles.
\newblock In {\em The Legacy of {John v. Neumann} ({Providence, 1990})},
  volume~50 of {\em Proc. Sympos. Pure Math.}, pages 151--178. American
  Mathematical Society, Providence, R.I., 1990.

\bibitem{Segal96b}
Irving~E. Segal.
\newblock Book review.
\newblock {\em Bull. Amer. Math. Soc. (N.S.)}, 33(4):459--465, 1996.
\newblock Review on book~\cite{Connes94}.

\bibitem{Segal96a}
Irving~E. Segal.
\newblock Rigorous covariant form of the correspondence principle.
\newblock In William~Arveson et~al., editor, {\em Quantization, Nonlinear
  Partial Differential Equations, and Operator Algebra ({Cambridge, MA},
  1994)}, volume~59 of {\em Proc. Sympos. Pure Math.}, pages 175--202. American
  Mathematical Society, Providence, R.I., 1996.

\bibitem{MTaylor86}
Michael~E. Taylor.
\newblock {\em Noncommutative Harmonic Analysis}, volume~22 of {\em Math. Surv.
  and Monographs}.
\newblock American Mathematical Society, Providence, R.I., {\noopsort{}}1986.

\end{thebibliography}
